\documentclass[10pt]{amsart}
\usepackage{amsmath, amsfonts, amsthm, amssymb, verbatim,dsfont}
\usepackage{graphicx}
\usepackage[mathcal]{euscript}
\usepackage{amstext}
\usepackage{dsfont} 
\usepackage{amssymb,mathrsfs}
\usepackage[usenames]{color}
\theoremstyle{plain}
        \newtheorem{theorem}{Theorem}[section]
        \newtheorem{proposition}[theorem]{Proposition}
        \newtheorem{lemma}[theorem]{Lemma}
        
\theoremstyle{definition}
        \newtheorem{definition}[theorem]{Definition}
        \newtheorem{remark}[theorem]{Remark}
\theoremstyle{plain}
                
\numberwithin{equation}{section}

\usepackage{amsmath}
\usepackage{verbatim}
\newcommand \be           {\begin{equation}}
\newcommand \ee            {\end{equation}}

\newcommand \RR           {\mathbb{R}}
\newcommand \NN           {\mathbb{N}}
\newcommand \ZZ           {\mathbb{Z}}
\newcommand \CC           {\mathbb{C}}

\newcommand \Pbold           {\mathbf{P}} 
\newcommand \PP \Pbold

\newcommand \del           \partial
\newcommand \eps            \epsilon
\newcommand \ubar       {{\overline u}}

\newcommand \loc        {{\mathrm{loc}}}

\DeclareMathOperator    \sgn {sgn}

\newcommand \ws     {\mathrel{\mathop{\rightharpoonup}\limits^{*}}}
\definecolor{gray}{gray}{0.4}
\newcommand{\unj}{u^n_j}
\newcommand{\unjj}{u^n_{j+1}}
\newcommand{\unjm}{u^n_{j-1}}
\newcommand{\unnj}{u^{n+1}_j}
\newcommand{\unnjj}{u^{n+1}_{j+1}}
\newcommand{\unnjm}{u^{n+1}_{j-1}}
\newcommand{\vnj}{v^n_j}
\newcommand{\vnjj}{v^n_{j+1}}
\newcommand{\vnjm}{v^n_{j-1}}
\newcommand{\vnnj}{v^{n+1}_j}
\newcommand{\vnnjj}{v^{n+1}_{j+1}}
\newcommand{\vnnjm}{v^{n+1}_{j-1}}

\newcommand{\vh}{v^h}
\newcommand{\uh}{u^h}
\newcommand{\vj}{v_j}
\newcommand{\vjm}{v_{j-1}}
\newcommand{\vjj}{v_{j+1}}
\newcommand{\uj}{u_j}
\newcommand{\ujm}{u_{j-1}}
\newcommand{\ujj}{u_{j+1}}

\newcommand{\ubj}{\ubar_j}

\newcommand{\ubjj}{\ubar_{j+1}}
\newcommand{\fjp}{f_{j,+}}
\newcommand{\fjm}{f_{j,-}}

\newcommand{\fjjm}{f_{j+1,-}}
\newcommand{\fjpm}{f_{j,\pm}}
\newcommand{\Gjp}{G_{j,+}}
\newcommand{\Gjm}{G_{j,-}}

\newcommand{\Gjjm}{G_{j+1,-}}
\newcommand{\Gjpm}{G_{j,\pm}}
\newcommand{\bGjp}{\bar{G}_{j,+}}
\newcommand{\bGjm}{\bar{G}_{j,-}}

\newcommand{\bGjjm}{\bar{G}_{j+1,-}}
\newcommand{\bGjpm}{\bar{G}_{j,\pm}}
\newcommand{\Hjp}{H_{j,+}}
\newcommand{\Hjm}{H_{j,-}}

\newcommand{\Hjpm}{H_{j,\pm}}

\newcommand{\etaj}{\eta_{j}}
\newcommand{\etajj}{\eta_{j+1}}
\newcommand{\fj}{f_{j}}
\newcommand{\fjj}{f_{j+1}}
\newcommand{\gj}{g'_{j}}
\newcommand{\gjj}{g'_{j+1}}

\DeclareMathOperator \supp {supp}
 
\DeclareMathOperator\Lip {Lip} 

\def\Xint#1{\mathchoice
{\XXint\displaystyle\textstyle{#1}}%
{\XXint\textstyle\scriptstyle{#1}}%
{\XXint\scriptstyle\scriptscriptstyle{#1}}%
{\XXint\scriptscriptstyle\scriptscriptstyle{#1}}%
\!\int}
\def\XXint#1#2#3{{\setbox0=\hbox{$#1{#2#3}{\int}$}
\vcenter{\hbox{$#2#3$}}\kern-.5\wd0}}

\def\medint{\Xint{-}}

\def\build#1_#2^#3{\mathrel{
\mathop{\kern 0pt#1}\limits_{#2}^{#3}}}

\begin{document}

\title[Numerical schemes for short wave long wave interaction equations]{Convergence of numerical schemes for short wave long wave interaction equations}

\author[P. Amorim \and M. Figueira]{Paulo Amorim$^1$ \and M\'ario Figueira$^1$}

\date{\today}

\footnotetext[1]{Centro de Matem\'atica e Aplica\c c\~oes
Fundamentais, FCUL, Av. Prof. Gama Pinto 2,
1649-003 Lisboa, Portugal. E-mail: {\tt pamorim@ptmat.fc.ul.pt, figueira@ptmat.fc.ul.pt}}

\begin{abstract} 
We consider the numerical approximation of a system of partial differential equations involving a 
nonlinear Schr\"odinger equation coupled with a hyperbolic conservation law. This system arises in 
models for the interaction of short and long waves. Using the compensated compactness method, 
we prove convergence of approximate solutions generated by semi-discrete finite volume type 
methods towards the unique entropy 
solution of the Cauchy problem. Some numerical examples are presented.
\end{abstract}

 
\maketitle

\section{Introduction}
\subsection{Interaction equations of short and long waves}
The nonlinear interaction between short waves and long waves has been studied in 
a variety of physical situations. In \cite{Benney}, D.J. Benney presents a general theory, deriving 
nonlinear differential systems involving both short and long waves. The short waves $u(x,t)$ are
described by a nonlinear Sch\-r\"o\-din\-ger equation and the long waves $v(x,t)$ satisfy a quasilinear 
wave equation, eventually with a dispersive term. 
In its most general form, the interaction is described by the nonlinear system
\be
\label{I.010} \left\{
\aligned 
& i \del_t u + i c_1 \del_x u + \del_{xx} u= \alpha u\,v + \gamma |u|^2 u
\\
& \del_t v +  c_2 \del_x v + \mu \del^3_x v + \nu \del_x v^2 = \beta \del_x (|u|^2),
\endaligned \right.
\ee
where $c_1,c_2, \alpha, \beta, \gamma, \mu$ and $\nu$ are real constants.

Systems of this type have attracted attention ever since the pioneering works of
M. Tsutsumi and S. Hatano \cite{TH1,TH2}, who set $c_2 = \mu = \nu = 0$ in \eqref{I.010} 
and proved the well-posedness of the Cauchy problem.
This result was later generalized by Bekiranov \emph{et al.} \cite{BOP}.
More recently, there has been renewed interest in systems of the type \eqref{I.010}. Here we just refer
to \cite{CorchoLinares} for the study of a Sch\-r\"o\-din\-ger--KdV system, 
and to \cite{DFF} for an in-depth study of several generalizations of \eqref{I.010}, among which is the
one concerning us in this work.

Regarding the numerical study of systems of type \eqref{I.010}, the only available works are, to the 
authors' knowledge, the papers \cite{AF} and \cite{Chineses}.
In \cite{AF}, the convergence (in $H^1\times L^2$) of a semi-discrete finite difference approximation 
of a version of the system \eqref{I.010} with $c_2 = \mu = \nu = 0$ is proved.

\subsubsection{A Sch\-r\"o\-din\-ger--conservation law system}
Here, following Dias \emph{et al.} \cite{DFF}, we consider the case where the long waves are 
modelled by a nonlinear hyperbolic conservation law. Thus in \cite{DFF} the following system is proposed,
\begin{subequations}
\begin{align}
\label{I.10}
&i\del_t u + \del_{xx} u =\, |u|^2 u +  \alpha g(v) u
\\
\label{I.20}
&\del_t v + \del_x f(v) = \alpha\del_x(g'(v) |u|^2)
\\
\label{I.25}
& u(x, 0) = u_0(x),\qquad  v(x, 0) =  v_0(x),
\end{align}
\end{subequations}
where $f \in C^2(\RR)$ and $g\in C^3(\RR)$ are real functions such that $f(0) = 0$ and $g'$ has 
compact support, and $\alpha>0$ is a constant. 
In that work, the authors establish the well-posedness of the Cauchy problem \eqref{I.10}--\eqref{I.25}, 
for an appropriate notion of entropy solution, provided $\alpha$ is small enough.

Comparing this model to the original formulation of Benney, we remark the appearance of a 
general interaction function $g(v)$ in the Sch\-r\"o\-din\-ger equation and of its derivative $g'(v)$ 
in the conservation law. This new function, it is argued in \cite{DFF}, gives a more physically realistic 
coupling between the two equations. The fact that $g'$ has compact support allows the authors 
of \cite{DFF}
to obtain a uniform bound for $\| v\|_\infty$ through a sort of maximum principle. This $L^\infty$
bound is fundamental in the analysis, and it is not clear whether a corresponding 
well-posedness result can be obtained for system \eqref{I.10}--\eqref{I.25} with $g(v) = v$, and for 
general $f$ (see, however \cite{DF} and \cite{ADFO} for some partial well-posedness and blow-up results, 
respectively).

Our objective here is to extend the numerical results in \cite{AF} to deal with the system 
\eqref{I.10}--\eqref{I.25}.
We will prove the convergence of a large class of schemes using
the compensated compactness method \cite{DiPerna,Tartar}, aided by energy inequalities 
which play the role 
of the classical entropy inequalities in the analysis of the method. 

We may summarize some of the main difficulties as follows. The equation \eqref{I.20} may be seen 
as a conservation law with a non-homogeneous flux-function given by $f(v) - 
g'(v) |u|^2$, or as a conservation law with the source term $\del_x(g'(v) |u|^2)$. However, to the 
best 
of our knowledge, the available convergence techniques for this kind of problem would require a 
pointwise 
bound on $|u|^2$ (see, for instance, \cite{BPV,CC}). For problem \eqref{I.10}--\eqref{I.25}, such a bound
is only obtained \emph{a posteriori}, and only as a consequence of the fact that $u\in H^1(\RR)$ and a 
Sobolev imbedding.

Therefore, the crucial bounds in Lemma \ref{UB-30} and Proposition \ref{UB-20} below
must employ new techniques which take into account the nonlinear coupling of the equations
\eqref{I.10},\eqref{I.20}, as was already observed in \cite{DF,DFF}. In the present paper we 
combine the techniques of \cite{DF,DFF} with those of the compensated compactness method
to derive new, stronger, compactness estimates which ensure the convergence of our 
numerical schemes.

\subsection{Outline of the paper}

In this paper, we establish the convergence of a class of numerical schemes to approximate the 
problem \eqref{I.10}--\eqref{I.25}. To this end, we use a semidiscrete
finite volume type scheme as an approximation of the quasilinear equation \eqref{I.20} 
and a semidiscrete finite-difference scheme for the nonlinear Sch\-r\"o\-din\-ger equation, 
\eqref{I.10}. The compensated compactness method allows us to prove
convergence towards the unique entropy solution, a result which had been announced (for
the particular case of the Lax-Friedrichs scheme) in \cite{AF}. 

The convergence result concerning problem \eqref{I.10}--\eqref{I.25} (Theorem \ref{20-10} below) relies on 
the compensated compactness method and on the crucial energy estimates in Lemma \ref{UB-30} 
and Proposition \ref{UB-20}. Interestingly, these estimates are finer than the corresponding ones in 
\cite{DFF}, and as a consequence our result does not require a smallness assumption on the coupling
parameter $\alpha$ in \eqref{I.10},\eqref{I.20}, and which we take equal to unity for simplicity. 
Also, no existence of solution is assumed \emph{a priori}, and thus our 
convergence proof is also a new existence proof for the problem \eqref{I.10}--\eqref{I.25}.


An outline of the paper follows. In Section \ref{Sec20} we define our numerical method and state
the main convergence result (Theorem \ref{20-10}). In Section~\ref{UB} we prove the main 
stability and energy estimates.
In Section~\ref{PC} we apply the compensated compactness method and the estimates of 
Section~\ref{UB} to deduce compactness of the approximate solutions. We then prove our 
convergence result. Finally, in Section~\ref{NE}, we present some numerical results
and discuss some open questions.

\section{The numerical scheme and statement of the main result}
\label{Sec20}

First of all, we recall from \cite{DFF} the notion of entropy solution to problem \eqref{I.10}--\eqref{I.25}.
\begin{definition}
\label{S-05}
Let $\eta(v)$ be a convex function (the \emph{entropy}), and define the \emph{entropy fluxes} $q_1,q_2$ by $q_1'(v) = \eta'(v) f'(v)$ and $q_2'(v) = \eta'(v) g''(v)$. We say that $(u,v) \in L^\infty_\loc( \RR \times [0,\infty))$ is an entropy solution to the problem \eqref{I.10}--\eqref{I.25} (with 
$\alpha = 1$) if for each entropy triplet $(\eta, q_1, q_2)$ we have:
\enumerate
\item $u \in L^\infty_\loc( [0,\infty) ; H^1(\RR)) \cap C([0, \infty) ; L^2(\RR))$, $u(0) = u_0$ in $L^2(\RR)$, and
\[
\aligned
& \iint_{\RR\times[0,\infty)} i u\, \del_t \theta + \del_x u\, \del_{x}\theta  + ( |u|^2 u + g(v) u) \theta \, dx dt =0
\endaligned
\]
for every $\theta \in C^\infty_0(\RR \times (0,\infty))$;
\item For every non-negative $\phi = \phi(x,t) \in C^\infty_0(\RR^2),$
\[
\aligned
& \iint_{\RR\times[0,\infty)}  \eta(v) \del_t \phi +  ( q_1(v) - q_2(v) |u|^2 ) \del_x \phi
\\ 
&\quad+ \big( \eta'(v) g'(v) - q_2(v) \big) \del_x |u|^2 \phi  \, dx dt + \int_\RR \eta(v_0(x)) \phi(x,0) \,dx \ge 0.
\endaligned
\]
\endenumerate
\end{definition}

\subsection{Finite volume schemes}

In this section, we define a class of semidiscrete finite volume schemes to approximate
\eqref{I.10}--\eqref{I.25}. 

Consider a uniform spatial grid $(x_j)_{j\in\ZZ}$ of element size $h = x_{j+1} - x_j$, where $h$ will be the discretization 
parameter. We seek functions $(\uh(t),\vh(t))$ where for each $t\ge0$ $\uh(t)$ and $\vh(t)$ are 
sequences $\uj(t) \in \CC$ and $\vj(t) \in \RR$.

For the nonlinear Sch\-r\"o\-din\-ger equation \eqref{I.10}, we use a standard finite difference scheme.
The term $\del_{xx} u$ is approximated by the usual discrete laplacian, denoted by
\[
(D^2_h \uh)_j \equiv D^2_h \uj = \frac{1}{h^2}( \ujj - 2\uj + \ujm).
\]
We will also need the difference quotient
\[
D_+ \uj = \frac{\ujj - \uj}{h}.
\]

For the second equation, \eqref{I.20}, we will consider a finite volume scheme, as follows. First, let us rewrite \eqref{I.20} in the form
\be
\label{15.10}
\del_t v + \del_x \big( f(v) - g'(v) |u|^2 \big) = 0.
\ee
Now for each $j\in \ZZ$ consider monotone, conservative, Lipschitz continuous numerical flux-functions
\[
\aligned
f_{j,+} (v_1, v_2), \quad f_{j,-} (v_1, v_2)  
\endaligned
\]
consistent with the function $f$, that is, $\fjp, \fjm$ verify for each $j\in\ZZ$,
$v_1, v_2, v \in \RR$, the conditions 
\begin{subequations}
\label{15.30}
\begin{align}
&\fjp (v_1,v_2) = - \fjjm (v_2,v_1)  && {\text{(Conservation)}}
\\
&\fjpm (v,v) = \pm f(v) 			&& \text{(Consistency)}
\\
& \del_1 \fjpm \ge 0, \quad \del_2 \fjpm \le 0  && \text{(Monotonicity),}
\end{align}
\end{subequations}
where $\del_1$ and $\del_2$ denote the derivative with respect to the first and second arguments,
respectively.

Next, consider another family $\bGjp(v_1, v_2)$, $\bGjm(v_1, v_2)$ of numerical flux-functions verifying the conditions 
in \eqref{15.30} but with $-g'(v)$ instead of $f(v)$. Thus, 
\begin{subequations}
\label{15.35}
\begin{align}
&\bGjp (v_1,v_2) = - \bGjjm (v_2,v_1) 
\\
&\bGjpm (v,v) = \mp g'(v) 
\\
& \del_1 \bGjpm \ge 0, \quad \del_2 \bGjpm \le 0.
\end{align}
\end{subequations}
Further, we consider a family of Lipschitz continuous functions
\[
\aligned
\Gjp( v_1, v_2 , a_1, a_2),
\endaligned
\]
satisfying the following conditions:
\begin{subequations}
\label{15.40}
\begin{align}
&\Gjp (v_1,v_2, a_1, a_2) = - \Gjjm (v_2,v_1, a_2, a_1) 
\\
&\Gjpm (v,v, a_1, a_2) = \mp  g'(v) \frac12(a_1 + a_2) 
\\
&\Gjpm (v_1, v_2 , a, a) = a\, \bGjpm(v_1, v_2) 
\\
& \del_1 \Gjpm \ge 0, \quad \del_2 \Gjpm \le 0.
\end{align}
\end{subequations}
Thus, if for instance $a = |u|^2$ for some complex number $u$, we find
\[
\Gjp (v,v,|u|^2, |u|^2) = - g'(v) |u|^2,
\]
so that this term gives a consistent approximation of the term $- g'(v) |u|^2.$

Finally, define $\Hjpm$ by
\be
\label{15.55}
\aligned
\Hjpm (v_1, v_2, a_1, a_2) = \fjpm (v_1, v_2) + \Gjpm(v_1, v_2, a_1, a_2).
\endaligned
\ee
In short, $\Hjpm$ are Lipschitz continuous, monotone, conservative, numerical flux functions, consistent 
with $f(v) - g'(v) |u|^2$ in the sense that
\be
\label{15.60}
\aligned
\\
\Hjpm (v, v, |u|^2, |u|^2) =\pm f(v) + \Gjpm(v,v,|u|^2, |u|^2) = \pm f(v) \mp g'(v) |u|^2.
\endaligned
\ee
For background on finite volume schemes, see \cite{EGH,GR}.

We now approximate the equation \eqref{15.10} by a semi-discrete finite volume method. We give
some initial data $u_{0j}, v_{0j}$ on $t = 0$; for instance, $u_{0j} = \frac1h\int_{x_j}^{x_{j+1}} u_0(x) dx$ 
(and similarly for $v_0$), so that the piecewise constant functions defined from $u_{0j}$, $v_{0j}$
converge strongly in $L^1_\loc(\RR)$ to $u_0, v_0$.  We then solve the infinite system of ODEs
\begin{subequations}
\label{15.50}
\begin{align}
\label{15.51}
& i \del_t \uj + D^2_h \uj = |\uj|^2 \uj + g(\vj) \uj,
\\
\label{15.52}
& \del_t v_j +\frac{1}{h} \big( \Hjp (\vj, \vjj, |\uj|^2, |\ujj|^2) + \Hjm (\vj, \vjm, |\uj|^2, |\ujm|^2) 
\big) = 0.
\end{align}
\end{subequations}

\subsection{Examples}
\subsubsection{A Lax-Friedrichs scheme}
Perhaps the simplest example of a numerical flux satisfying the above conditions is given
by the following Lax--Friedrichs type scheme. Let $\lambda, \gamma$ be some constants and set
\be
\label{LF-01}
\aligned
\fjpm(v_1,v_2) &\equiv \mathbf{f}_\pm (v_1,v_2) = \pm\frac12 \big( f(v_1) + f(v_2) \big) 
\mp \frac1{2\lambda} (v_2 - v_1),
\\
\Gjpm( v_1,v_2, a_1, a_2) &\equiv \mathbf{G}_\pm (v_1,v_2,a_1, a_2) 
\\
&\qquad =
\pm\frac12 \big( -g'(v_1) a_1 - g'(v_2) a_2 \big) 
\mp \frac1{2\gamma} (v_2 - v_1) \frac12 (a_1 + a_2).
\endaligned
\ee
It is easy to check that, under an appropriate CFL-type condition on $\lambda$ and $\gamma$, 
the numerical fluxes $\mathbf{f}_\pm$ and $\mathbf{G}_\pm$ verify the assumptions of 
the previous section.

\subsubsection{A Godunov scheme}
Here we consider the following straightforward generalization of the classical Godunov scheme. Setting
\[
\aligned
H_{j,\pm} \equiv\mathbf{H}_{\pm} (v_1, v_2, a_1, a_2) =
\left\{ 
\aligned
&\min_{v_1 \le s \le v_2} \pm \big( f(s) - g'(s) \frac12(a_1+a_2) \big), && v_1\le v_2,
\\
&\max_{v_2 \le s \le v_1} \pm \big( f(s) - g'(s) \frac12(a_1+a_2) \big), && v_2\le v_1,
\endaligned\right.
\endaligned
\]
one can then verify that the consistency, conservation and monotonicity conditions above are met.
\subsection{Convergence result}

We are now ready to state the main result of this paper, whose proof 
will be the object of Section~\ref{PC}.

In what follows, we suppose that $f$ and $g$ verify a standard non-degeneracy condition:
\[
\forall \, \kappa >0, \quad\text{the set} \quad \{ s : f''(s) - \kappa g'''(s) \neq 0 \} \quad
\text{is dense in $\RR$}.
\]
This condition is needed in order to apply the compensated compactness method.

\begin{theorem}
\label{20-10}
Let $(\uh,\vh)$ be defined by a semidiscrete scheme \eqref{15.50}, verifying all the 
assumptions \eqref{15.30}--\eqref{15.40}.  Then
there exist functions 
\[
u \in C([0,\infty); H^1(\RR)), \qquad v\in L^\infty(\RR\times [0,\infty)),
\] 
solutions of the Cauchy problem \eqref{I.10}--\eqref{I.25} (with $\alpha =1$) in the sense of 
Definition \ref{S-05} such that, up to a subsequence, $(\uh,\vh)$ converge to 
$(u,v)$ in $L^1_\loc(\RR\times [0,\infty))$.
\end{theorem}

\section{Uniform bounds and energy estimates}
\label{UB}
In this section we prove the crucial uniform estimates of Propositions~\ref{UB-10}, \ref{UB-20}, and 
Lemma~\ref{UB-30}. We begin with a key property of finite volume schemes.

Let $\fjp$ be a numerical flux-function consistent with a function $f$, let $v_1, v_2 \in\RR$, and
let $\eta$ be a smooth convex function. Define 
the \textbf{viscosity} as the quantity
\be
\label{Q1}
\int_{v_1}^{v_2} \eta''(v)( f(v) - \fjp(v_1, v_2)) dv.
\ee
The following property, although classical, is very important in the analysis of the scheme 
(cf. \cite{Tadmor,GR}). It is an easy consequence of the monotonicity and consistency properties of 
the flux functions, and so we omit its proof.
\begin{lemma}
\label{L-10}
For each $v_1,v_2 \in \RR$ and each smooth convex function, the viscosity 
\eqref{Q1} is everywhere non-negative.
In particular, we have
\[
\int_{v_1}^{v_2} \eta''(v)( f(v) - \fjp(v_1, v_2)) dv \ge 0,
\]
\[
\int_{v_1}^{v_2} \eta''(v)( -g'(v) - \bGjp(v_1, v_2)) dv \ge 0
\]
for every $v_1, v_2 \in \RR$.
\end{lemma}

In what follows, if (say) $\vh$ is a sequence, we will not distinguish between the discrete norm 
of $\vh$ in $l^p(\ZZ)$ spaces its and continuous norm $\|\vh\|_p$ defined as the usual 
$L^p$ norm of the piecewise constant function $\vh(x) = \vj$ if $x\in [ x_j, x_{j+1}]$. Thus if
$p \in [1, \infty)$,
\[
\|\vh\|_p = \Big( \sum_{j\in\ZZ} h|\vj|^p\Big)^{1/p},
\]
with the usual modification if $p=\infty$.

We now prove an $L^\infty$ bound for $\vh$ which is essential in all our analysis, as well as 
the conservation of the $L^2$ norm of $\uh$.
\begin{proposition}
\label{UB-10}
Let $(\uh, \vh)$ be defined by the finite volume method \eqref{15.50}, for some initial data
$u_{0j} \in l^2(\ZZ)$, $v_{0j} \in l^2(\ZZ)$. Then, for each $h>0$, 
there is a unique global solution of \eqref{15.50}. Moreover, for some $M = M(v_0, g) >0$
independent of $h$ (where $g$ is the coupling function in \eqref{I.10}), we have
\be
\label{UB.10}
\frac{d}{dt}\| \uh(t) \|_2^2  =0
\ee
and
\be
\label{UB.20}
\| \vh(t) \|_\infty \le M.
\ee
\end{proposition}
\begin{proof}
First, the local in time existence of solution of \eqref{15.50} is ensured by a standard 
fixed point argument in the discrete space $l^2(\ZZ)$. The argument is very similar to the one in
\cite{DFF}, and so we omit the 
details.
 
Next, we sketch the proof of estimate \eqref{UB.10}, which consists of a straightforward calculation. 
From the scheme \eqref{15.51}, multiply
by $h\uj$, sum in $j \in \ZZ$ and sum by parts the term with $D^2_h u_j$ 
to obtain that
\[ i \del_t \sum_{j\in\ZZ}h | u_j|^2 \in \RR, \]
which immediately gives \eqref{UB.10}.

Since, for fixed $h>0$, one has $l^2(\ZZ) \subset l^\infty(\ZZ)$, the local solution will be 
global in time if we prove the bound \eqref{UB.20}.

Following \cite{DFF}, consider 
$M'$ such that $\supp g' \subset [-M', M']$. Then, setting $M= \max\{ M', \| \vh(0)\|_\infty\}$,
we will see that $\|\vh \|_\infty \le M$. First, we consider a perturbed problem which coincides
with \eqref{15.50}, except that the term $-\eps \sgn \vh$ is added to the right-hand side of the 
second equation, for some small $\eps>0$. Let $(u^{h,\eps}, v^{h,\eps})$ be the solution 
of this problem, whose existence is obtained in the same way as for the unperturbed problem.

In what follows, let us consider only the terms 
involving the coupling function $g'$, since the treatment of the terms involving the
flux $f$ is standard and can be found for instance in \cite{GR}.
Thus, suppose that there is a first $j\in\ZZ$ and a first 
instant $t>0$ such that $\vj^\eps(t) =M$ (the case $\vj^\eps(t) = -M$ is similar). Then the assumptions
on the coupling function $g$, namely the compact support of $g'$, imply that the term with 
$\Gjpm$ appearing in the definition of the scheme \eqref{15.52} may be written as
\[
\aligned
& - \frac{1}{h} \big( \Gjp (\vj^\eps, \vjj^\eps, |\uj^\eps|^2, |\ujj^\eps|^2) + \Gjm (\vj^\eps, \vjm^\eps, |\uj^\eps|^2, |\ujm^\eps|^2) 
\big)
\\
&\quad = - \frac{1}{h} \big( \Gjp (\vj^\eps, \vjj^\eps, |\uj^\eps|^2, |\ujj^\eps|^2) +  \Gjp (\vj^\eps, \vj^\eps, |\uj^\eps|^2, |\ujj^\eps|^2) 
\big)
\\
& \qquad -\frac{1}{h} \big( \Gjm (\vj^\eps, \vjm^\eps, |\uj^\eps|^2, |\ujm^\eps|^2) + \Gjm (\vj^\eps, \vj^\eps, |\uj^\eps|^2, |\ujm^\eps|^2) 
\big),
\endaligned
\]
due to $g'(\vj^\eps) = 0$ and the consistency conditions in \eqref{15.35},\eqref{15.40}. Now divide and
multiply the first sum by $\vjj^\eps - \vj^\eps$ (which is non-positive since $v^{h,\eps}$ 
attains its maximum value at 
$j$), and the second sum by $\vjm^\eps - \vj^\eps$. 
The monotonicity condition in \eqref{15.40} then ensures 
that the resulting terms have the right sign so that (along with the treatment of the terms involving
$\fjp$, which we omit)
\[
\del_t \vj^\eps(t)  \le -\eps <0,
\]
which is absurd, since we assumed that $v^{h,\eps}$ attained a maximum value at $j, t$. 
Thus, $|v^{h,\eps}|<M$.

It is now necessary to take the limit $\eps \to0$ (with $h$ fixed) to deduce the desired estimate
\eqref{UB.20}. We 
obtain from the scheme \eqref{15.52} using the regularity of the numerical flux functions $\Hjpm$,
\be
\label{max1}
\aligned
\sup_{j\in\ZZ} |v_j^\eps (t) - \vj(t) | \le C(h) \int_0^t \sup_{j\in\ZZ} |v_j^\eps  - \vj|
+ \sup_{j\in\ZZ} \big| |u_j^\eps|^2 - |\uj|^2 \big| \,ds + \eps t,
\endaligned
\ee
where $C(h)$ is some (unbounded) function of $h$ which may change from line to line. 

Note that the estimate \eqref{UB.10} remains valid with $\eps >0$ with 
exactly the same proof, that is, we have $\| u^{h,\eps}(t) \|_2 =\| u^{h}(t) \|_2 =C.$ 
From
$l^2 \subset l^\infty$ we get immediately $\| u^{h,\eps}\|_\infty \le h^{-1/2}\| u^{h,\eps}\|_2$.
Thus, from \eqref{max1},
\[
\aligned
\sup_j \big| |u_j^\eps|^2 - |\uj|^2 \big| &\le 2\sup_j| \uj^\eps| | \uj^\eps - \uj|
\\
& \le 2 h^{-1/2} \| \uh\|_2 \sup_j |\uj^\eps - \uj|
\\
&\le C(h) \sup_j |u_j^\eps - \uj|,
\endaligned
\]
and using the scheme \eqref{15.51} gives in turn
\[
\aligned
\sup_{j\in\ZZ} |u_j^\eps - \uj| \le C(h) \int_0^t \sup_{j\in\ZZ} |u_j^\eps - \uj| 
+ \sup_{j\in\ZZ} | v_j^\eps - \vj| \,ds.
\endaligned
\]
Gronwall's lemma yields
\[
\aligned
\sup_{j\in\ZZ} |u_j^{\eps} - \uj| \le C(h,t)\sup_{j\in\ZZ} | v_j^{\eps} - \vj|,
\endaligned
\]
which we substitute in \eqref{max1} to get (after applying Gronwall again)
\[
\aligned
\sup_{j\in\ZZ} |v_j^{\eps} (t) - \vj(t) | \le C(h,t) \eps 
\endaligned
\]
which goes to zero with $\eps$, for $h$ and $t$ fixed. Therefore we deduce \eqref{UB.20}.
This completes the proof of Proposition~\ref{UB-10}.
\endproof

We now prove an energy inequality for the second equation of \eqref{15.50}, involving 
the viscosity \eqref{Q1}, which will
be essential to establish Proposition~\ref{UB-20}.

\begin{lemma}
\label{UB-30}
Let $(\uh, \vh)$ be defined by the scheme \eqref{15.50} and 
let $\eta$ be a smooth convex function. Then, the 
estimate
\be
\label{UB.100}
\aligned
\| \eta(\vh(t)) \|_1  {}&+ \int_0^t \sum_{j\in\ZZ}\int_{\vj}^{\vjj} \eta''(v)( f(v) - \fjp(\vj, \vjj)) dv 
\\
&\qquad+ |\uj|^2 \int_{\vj}^{\vjj} \eta''(v)( -g'(v) - \bGjp(\vj, \vjj)) dv  \, ds
\\
&\le c +c \int_0^t \| D_+ \uh(s)\|_2\, ds.
\endaligned
\ee
holds for some constant $c >0$ independent of $h$.
\end{lemma}
\begin{proof}
The scheme \eqref{15.50} reads
\[
\aligned
\del_t v_j &+ \frac{1}{h} \big( \fjp (\vj, \vjj) + \fjm (\vj, \vjm) 
\big) 
\\
& + \frac{1}{h} \big( \Gjp (\vj, \vjj, |\uj|^2, |\ujj|^2) + \Gjm (\vj, \vjm, |\uj|^2, |\ujm|^2) 
\big)
= 0.
\endaligned
\]
Now multiply by $h\eta'(\vj)$ and sum over $j \in \ZZ$ to obtain, with obvious notation,
\[
\aligned
& \del_t\sum_{j\in\ZZ}h \eta(\vj) + \sum_{j\in\ZZ} \eta'(\vj) \big( \fjp (\vj, \vjj) + \fjm (\vj, \vjm) 
\big) 
\\
& + \sum_{j\in\ZZ} \eta'(\vj) \big( \Gjp (\vj, \vjj, |\uj|^2, |\ujj|^2) + \Gjm (\vj, \vjm, |\uj|^2, |\ujm|^2) 
\big)
\\
& \equiv  \del_t\sum_{j\in\ZZ}h \eta(\vj) + A + B = 0.
\endaligned
\]
First, summation by parts gives
\[
A = - \sum_{j\in\ZZ} (\eta'(\vjj) - \eta'(\vj)) \fjp(\vj, \vjj) = -\sum_{j\in\ZZ} \int_{\vj}^{\vjj}  
\eta''(v)\fjp(\vj, \vjj) \,dv.
\]
Next, let $q_1(v)$ be such that $f'\eta' = q_1'$. We have
\[
\aligned
0 &= \sum_{j\in\ZZ} q_1(\vjj) - q_1(\vj) = \sum_{j\in\ZZ} \int_{\vj}^{\vjj} \eta'(v) f'(v) \, dv 
= -\sum_{j\in\ZZ} \int_{\vj}^{\vjj} \eta''(v) f(v) \, dv,
\endaligned
\]
and therefore
\be
\label{AA}
\aligned
A &=  - \sum_{j\in\ZZ} \int_{\vj}^{\vjj}  \eta''(v)(\fjp(\vj, \vjj) -  f(v)) \, dv .
\endaligned
\ee
Now consider the term $B$. First, we find
\[
\aligned
B &= - \sum_{j\in\ZZ} (\eta'(\vjj) - \eta'(\vj)) \Gjp(\vj, \vjj, |\uj|^2, |\ujj|^2)
\\
&= - \sum_{j\in\ZZ} \int^{\vjj}_{\vj} \eta''(v) \Gjp(\vj, \vjj, |\uj|^2, |\ujj|^2) \, dv.
\endaligned
\]
Next, setting $q_2' = \eta' g''$ gives
\[
\aligned
0 &= \sum_{j\in\ZZ} q_2(\vjj) |\ujj|^2 - q_2(\vj) |\uj|^2 
\\
&= \sum_{j\in\ZZ} (q_2(\vjj) - q_2(\vj) )|\uj|^2 + 
q_2(\vjj) (|\ujj|^2 - |\uj|^2 )
\\
& = \sum_{j\in\ZZ} |\uj|^2 \int_{\vj}^{\vjj} g''(v) \eta'(v) \, dv + q_2(\vjj) (|\ujj|^2 - |\uj|^2),
\endaligned
\]
and thus, integrating by parts,
\[
\aligned
 0& = -\sum_{j\in\ZZ} |\uj|^2 \int_{\vj}^{\vjj} g'(v) \eta''(v) \, dv +\sum_{j\in\ZZ} |\uj|^2 
\big( g'(\vjj) \eta'(\vjj) - g'(\vj) \eta'(\vj) \big)
\\
&\qquad + \sum_{j\in\ZZ} q_2(\vjj) (|\ujj|^2 - |\uj|^2)
\\
& =  -\sum_{j\in\ZZ} |\uj|^2 \int_{\vj}^{\vjj} g'(v) \eta''(v) \, dv -
\big( q_2(\vjj) - g'(\vjj)\eta'(\vjj) \big) (|\ujj|^2 - |\uj|^2).
\endaligned
\]
This gives
\[
\aligned
B &= - \sum_{j\in\ZZ} \int_{\vj}^{\vjj} \eta''(v) \big( \Gjp(\vj, \vjj, |\uj|^2, |\ujj|^2)
+ |\uj|^2  g'(v) \big) \,dv
\\
&\qquad + \sum_{j\in\ZZ} \big( q_2(\vjj) - g'(\vjj)\eta'(\vjj) \big) (|\ujj|^2 - |\uj|^2)
\\
& =  - \sum_{j\in\ZZ} |\uj|^2 \int_{\vj}^{\vjj}\eta''(v) \big( \bGjp(\vj, \vjj)
+g'(v) \big) \, dv 
\\
&\qquad
+ \sum_{j\in\ZZ} \big( q_2(\vjj) - g'(\vjj)\eta'(\vjj) \big) (|\ujj|^2 - |\uj|^2)
\\
&\qquad +  \sum_{j\in\ZZ} \int_{\vj}^{\vjj} \eta''(v) \,dv \Big( |\uj|^2  \bGjp(\vj, \vjj) 
- \Gjp(\vj, \vjj, |\uj|^2, |\ujj|^2) \Big).
\endaligned
\]
The first sum in the right-hand side of the previous equality is the sum appearing in \eqref{UB.100}.
The other two sums are bounded by
\[
\aligned
C \sum_{j\in\ZZ} \big| |\uj|^2 - |\ujj|^2 \big|  
\endaligned
\]
by appealing to $|\uj|^2  \bGjp(\vj, \vjj) = \Gjp(\vj, \vjj, |\uj|^2, |\uj|^2)$ (cf. \eqref{15.40}), the regularity of the numerical flux functions $\Gjp$, the $L^\infty$ bound on $\vj$, \eqref{UB.20}, 
and the smoothness of $g$ and $\eta$. Finally, in view of the bound on the $L^2$ norm of $\uh$, 
\eqref{UB.10}, we find
\be
\label{UB.101}
\aligned
\sum_{j\in\ZZ} \big| |\uj|^2 - |\ujj|^2 \big|  &\le \sum_{j\in\ZZ} |\ujj - \uj| | \ubjj + \ubj| 
\\
&\le 2 \Big( \sum_{j\in\ZZ} \frac1h| \ujj -\uj|^2 \sum_{j\in\ZZ} h|\uj|^2 \Big)^{1/2}
\\
&= 2 \| D_+ \uh \|_2 \| \uh\|_2 = c  \| D_+ \uh \|_2.
\endaligned
\ee
Thus 
\[
B =  \sum_{j\in\ZZ} |\uj|^2 \int_{\vj}^{\vjj}\eta''(v) \big( -\bGjp(\vj, \vjj)
-g'(v) \big) \, dv  + B_1,
\]
with $|B_1| \le c \| D_+ \uh\|_2$. This estimate and \eqref{AA} give
\eqref{UB.100}, after integration on $(0,t)$. This completes
the proof of Lemma~\ref{UB-30}.
\end{proof}

The following is the crucial \emph{a priori} estimate in our analysis.
\begin{proposition}
\label{UB-20}
Let $(\uh, \vh)$ be defined by the scheme \eqref{15.50}. Then,
there exist  
non-negative functions $a(t),$ $b(t)$ continuous on $[0, \infty)$ such that
\be
\label{UB.25}
\| \vh (t) \|_2 \le b(t),
\ee
\be
\label{UB.27}
 \| D_+ \uh(t) \|_2  \le a(t).
\ee
\end{proposition}
\begin{remark}
In view of the estimate \eqref{UB.27}, the viscosity estimate \eqref{UB.100} implies
\be
\label{UB.28}
\aligned
\int_0^t & \sum_{j\in\ZZ}\int_{\vj}^{\vjj} \eta''(v)( f(v) - \fjp(\vj, \vjj)) dv \, ds
\\
&\qquad+ \int_0^t |\uj|^2 \int_{\vj}^{\vjj} \eta''(v)( -g'(v) - \bGjp(\vj, \vjj)) dv  \, ds
\le c(t)
\endaligned
\ee
for any smooth convex function $\eta$ and some locally bounded function $c(t)$.
\end{remark}
\begin{remark}
In the case of the Lax--Friedrichs scheme \eqref{LF-01}, the estimate \eqref{UB.28} actually implies
an explicit bound on the ``quadratic total variation" of $\vh$, as was observed in \cite{AF}. 
More precisely, we have for some $k>0$
\be
\label{LF-2}
\aligned
\| \vh(t)\|_2^2 + k  \int_0^t
\sum_{j\in\ZZ}(1+ |\uj|^2) (\vjj - \vj)^2   \,ds\le   c(t).
\endaligned
\ee
This can be seen by observing that for the Lax--Friedrichs scheme, one has 
(for $v_1 \le v \le v_2$; the other case is similar)
\[
\aligned
\frac{ f(v) - \mathbf{f}_+(v_1,v_2)}{v_2 - v_1}  &= \frac{ \mathbf{f}_+(v,v) - 
\mathbf{f}_+(v_1,v_2)}{v_2 - v_1} 
\\
&\ge  \frac{ \mathbf{f}_+(v_1,v) - \mathbf{f}_+(v_1,v_2)}{v_2 - v_1} + k
\\
&\ge  \frac{ \mathbf{f}_+(v_1,v_2) - \mathbf{f}_+(v_1,v_2)}{v_2 - v_1} + k = k >0,
\endaligned
\]
for some $k>0$. This follows from the monotonicity of the scheme, and $k$ appears due to the 
particular form of the Lax--Friedrichs scheme and a CFL condition $\lambda |f'| < 1$. Thus 
(see \eqref{UB.28}, with $\eta'' =1$)
\be
\label{LF-3}
\aligned
\int_{v_1}^{v_2}  f(v) - \mathbf{f}_+(v_2, v_1) dv \ge \int_{v_1}^{v_2} k(v_2 - v_1)\, dv
= k (v_2 - v_1)^2.
\endaligned
\ee
This, together with a similar estimate for the numerical flux $\mathbf{G}_+$ in \eqref{LF-01}
gives \eqref{LF-2}. In the general case, however, the estimate \eqref{LF-3} is 
only valid for $k=0$ (for example, the Godunov scheme in the neighborhood of
a sonic point fails to verify \eqref{LF-3} with $k>0$). Thus, in general there is no bound on the
quadratic total variation as in \eqref{LF-2}, and the appropriate estimate is
\eqref{UB.28}.
\end{remark}

\proof[Proof of Proposition \ref{UB-20}]
First, of all, recall the 
Gagliardo--Niremberg--Sobolev inequalities which we will use throughout. If $\phi \in l^2(\ZZ)$, then
\begin{align}
\label{M.10}
&\| \phi \|_\infty \leq C \| \phi\|_{2}^{1/2} \|D_+\phi\|_{2}^{1/2}    \\
\label{M.20}
& \| \phi \|_{4} \leq C \| \phi\|_{2}^{3/4} \|D_+\phi\|_{2}^{1/4}.
\end{align}
Now take \eqref{15.51}, multiply
by $h\del_t\bar u_j$ and sum over $j \in \ZZ$ to obtain, using straightforward computations very
similar to the ones detailed in \cite[Lemma 2.2]{DFF},
\be
\label{UB.130}
\aligned
 \frac12 \| D_+ \uh\|_2^2 + \frac14 \|\uh \|_4^4 + \frac12 \sum_{j\in\ZZ}
h g(\vj) |\uj|^2 
 = c + \int_0^t \sum_{j\in\ZZ}h g'(\vj) \del_t \vj |\uj|^2 ds.
\endaligned
\ee
Let us now consider the term on the right-hand side of \eqref{UB.130}.
Using the second equation of the scheme \eqref{15.52}, and the conservation properties of the 
numerical flux functions, 
we find (omitting the arguments of the numerical flux functions)
\be
\label{UB.132}
\aligned
&\sum_{j\in\ZZ}h g'(\vj) \del_t \vj |\uj|^2  = - \sum_{j\in\ZZ} g'(\vj) |\uj|^2
\big( \fjp + \fjm + \Gjp + \Gjm \big)
\\
& \quad= \sum_{j\in\ZZ} \big( g'(\vjj) |\ujj|^2 - g'(\vj) |\uj|^2 \big) (\fjp + \bGjp).
\endaligned
\ee
First, we have
\be
\label{UB.1320}
\aligned
\sum_{j\in\ZZ}  \big( g'(\vjj) |\ujj|^2 - g'(\vj) |\uj|^2 \big) \fjp  &
 = \sum_{j\in\ZZ} |\uj|^2 ( g'(\vjj) - g'(\vj) )  \fjp 
\\& \hspace{-20mm}+ \sum_{j\in\ZZ}  g'(\vjj)( |\ujj|^2 -  |\uj|^2 ) \fjp =: A_1 + A_2.
\endaligned
\ee
We will write $A_1$ in terms of the viscosity in order to use the estimate \eqref{UB.28}.
\[
\aligned
A_1 &= \sum_{j\in\ZZ} |\uj|^2 \int_{\vj}^{\vjj} g''(v) \fjp \,dv
\\
& = \sum_{j\in\ZZ} |\uj|^2 \int_{\vj}^{\vjj} g''(v) (\fjp - f(v)) \,dv
+\sum_{j\in\ZZ} |\uj|^2 \int_{\vj}^{\vjj} g''(v)f(v) \,dv.
\endaligned
\]
It is an easy consequence of the monotonicity and consistency properties of the scheme 
that for $s$ between $\vj$ and $\vjj$, the quantity $\fjp(\vj, \vjj) - f(s)$ does not 
change sign (cf. Lemma \ref{L-10}). Therefore, we may write for some intermediate
value $c$ (observing that from Lemma~\ref{L-10} the 
integral below has the right sign)
\[
\aligned
\Big|\int_{\vj}^{\vjj} g''(v) (\fjp - f(v)) \,dv \Big| = |g''(c)| \int_{\vj}^{\vjj} (f(v) - \fjp) \,dv.
\endaligned
\]
Moreover, defining a function $\mathcal{F}$ by $\mathcal{F}'(v) = g''(v) f(v)$, we find
\[
\aligned
\sum_{j\in\ZZ} |\uj|^2 \int_{\vj}^{\vjj} g''(v)f(v)  \,dv 
&= \sum_{j\in\ZZ} |\uj|^2 (\mathcal{F}(\vjj) - \mathcal{F}(\vj) )  \,dv
\\
& = - \sum_{j\in\ZZ} \big( |\ujj|^2 - |\uj|^2 \big) \mathcal{F}(\vjj)   \,dv, 
\endaligned
\]
and thus, by the uniform bound for $\vh$ \eqref{UB.20} and the Sobolev inequality \eqref{M.10},
\[
\aligned
A_1 &\le  \sum_{j\in\ZZ} |\uj|^2 \sup|g''| \int_{\vj}^{\vjj}  (f(v) - \fjp ) \,dv
+\sum_{j\in\ZZ}  \sup|\mathcal{F}|  \big| |\ujj|^2 - |\uj|^2 \big|
\\
&\le c \|D_+ \uh\|_2 \sum_{j\in\ZZ} \int_{\vj}^{\vjj}  (f(v) - \fjp ) \,dv + c\| D_+ \uh\|_2.
\endaligned
\]
We may now apply Lemma \ref{UB-30} with $\eta(v) = v^2/2$ to get after integration on $(0,t)$
\[
\aligned
\int_0^t A_1(s) \, ds & \le \sup_{(0,t)} \|D_+ \uh(s)\|_2 \Big( c + c \int_0^t \| D_+ \uh(s)\|_2 ds
\Big)\\
&\quad+ c \int_0^t \| D_+ \uh(s) \|_2 ds.
\endaligned
\]
The term $A_2$ in \eqref{UB.1320} is estimated more simply using the uniform bound for 
$\vh$ \eqref{UB.20} to give 
after integration on $(0,t)$
\[
\aligned
\int_0^t A_2(s) \,ds &\le c \int_0^t \|D_+ \uh(s) \|_2 ds.
\endaligned
\]
Thus, \eqref{UB.1320} becomes
\be
\label{UB.133}
\aligned
\int_0^t \sum_{j\in\ZZ} & \big( g'(\vjj) |\ujj|^2 - g'(\vj) |\uj|^2 \big) \fjp  ds\\
&\le \sup_{(0,t)} \|D_+ \uh\|_2 \Big( c + c \int_0^t \| D_+ \uh\|_2 ds
\Big)+ c \int_0^t \| D_+ \uh \|_2 ds.
\endaligned
\ee
The remaining term in \eqref{UB.132} gives
\be
\label{UB.140}
\aligned
\sum_{j\in\ZZ}  \big( g'(\vjj) |\ujj|^2 - g'(\vj) |\uj|^2 \big) \Gjp  &
 = \sum_{j\in\ZZ} |\uj|^2 ( g'(\vjj) - g'(\vj) )  \Gjp 
\\& \hspace{-20mm}+ \sum_{j\in\ZZ}  g'(\vjj)( |\ujj|^2 -  |\uj|^2 ) \Gjp =: B_1 + B_2.
\endaligned
\ee
We have
\be
\label{UB.150}
\aligned
B_1 &= \sum_{j\in\ZZ} |\uj|^4 ( g'(\vjj) - g'(\vj) )  \bGjp 
\\
& \quad+ \sum_{j\in\ZZ} |\uj|^2 ( g'(\vjj) - g'(\vj) ) \big( \Gjp - |\uj|^2 \bGjp\big) =:
B_{11} + B_{12}.
\endaligned
\ee
Again, we use estimate \eqref{UB.28}, and so we need the viscosity to appear. We find
\[
\aligned
B_{11} &= \sum_{j\in\ZZ} |\uj|^4 \int_{\vj}^{\vjj} g''(v)  \big( \bGjp + g'(v) \big)dv
 - \sum_{j\in\ZZ}  |\uj|^4 \int_{\vj}^{\vjj} g''(v) g'(v) dv
\\
& \le \| \uh\|_\infty^2 \sum_{j\in\ZZ} \sup|g''| |\uj|^2 \int_{\vj}^{\vjj} (-\bGjp - g'(v)) dv
\\
& \quad+ \sum_{j\in\ZZ} \frac12 (g'(\vjj))^2 \big| |\ujj|^4 - |\uj|^4 \big|
\endaligned
\]
which gives, using Lemma \ref{UB-30} and the Sobolev inequality \eqref{M.10}
\[
\aligned
\int_0^t B_{11} \,ds & \le c\sup_{(0,t)} 
\| D_+ \uh \|_2 \int_0^t \sum_{j\in\ZZ} |\uj|^2 \int_{\vj}^{\vjj} (-\bGjp - g'(v)) dv ds
\\
& \quad
 + c\int_0^t \| \uh\|_\infty^2 \sum_{j\in\ZZ} \big| |\ujj|^2 - |\uj|^2 \big| ds
\\
& \le \sup_{(0,t)} \|D_+ \uh\|_2 \Big( c + c \int_0^t \| D_+ \uh\|_2 ds
\Big).
\endaligned
\]
The term $\int_0^t B_{12} \,ds$ from \eqref{UB.150} yields using \eqref{UB.20} and \eqref{M.10} (recall
the relation between $\Gjp$ and $|\uj|^2 \bGjp$ in \eqref{15.40})
\[
\aligned
\int_0^t B_{12}\, ds \le \int_0^t \sum_{j\in\ZZ}|\uj|^2 \Lip G \big| |\ujj|^2 - |\uj|^2 \big|\,ds
\le c \sup_{(0,t)} \|D_+ \uh\|_2  \int_0^t \| D_+ \uh\|_2 ds.
\endaligned
\]
As for $B_2$ in \eqref{UB.140}, we find immediately using the uniform bound for $\vh$ \eqref{UB.20}
\[
\aligned
\int_0^t B_{2}\, ds \le c \int_0^t \|D_+ \uh(s) \|_2 ds.
\endaligned
\]
Thus, \eqref{UB.132} becomes
\be
\label{UB.160}
\aligned
&\int_0^t \sum_{j\in\ZZ}h g'(\vj) \del_t \vj |\uj|^2 ds \\
& \quad\le \sup_{(0,t)} \|D_+ \uh\|_2 \Big( c + c \int_0^t \| D_+ \uh\|_2 ds
\Big)+ c \int_0^t \| D_+ \uh \|_2 ds,
\endaligned
\ee
and, after observing that estimates \eqref{UB.10} and \eqref{UB.20} imply
\[
 \sum_{j\in\ZZ} h |g(\vj)| |\uj|^2  \le c,
\]
we conclude from \eqref{UB.130}, \eqref{UB.133} and \eqref{UB.160} that
\[
\aligned
\frac12 \| D_+ \uh(t)\|_2^2 + \frac14 \|\uh(t) \|_4^4
& \le c  + c \int_0^t \|D_+ \uh \|_2 ds 
\\
&\qquad+ \sup_{(0,t)} \|D_+ \uh \|_2 \Big( c +c \int_0^t \|D_+ \uh \|_2 ds\Big).
\endaligned
\]
In particular, setting $\beta(t) := \sup_{(0,t)}\big( 1 + \| D_+ \uh(t)\|_2^2 \big)$ this
easily implies
\[
\aligned
\beta(t)
& \le c + c \int_0^t \beta(s)^{1/2} ds 
\\
&\qquad+ \beta(t)^{1/2} \Big( c +c \int_0^t \beta(s)^{1/2} ds\Big).
\endaligned
\]
The estimate \eqref{UB.27} follows by Gronwall's lemma applied to the function
$\beta(t)^{1/2}$. To prove \eqref{UB.25}, simply apply \eqref{UB.27} to the estimate
\eqref{UB.100} in Lemma~\ref{UB-30}, with $\eta(v) = v^2/2$. This completes the proof of Proposition~\ref{UB-10}.
\end{proof}

\section{Proof of convergence}
\label{PC}
In this section we prove Theorem \ref{20-10},
relying on the compensated compactness method \cite{Murat,Tartar}. 
According to this method, the strong compactness of a sequence of 
approximate solutions $(v^h)$ is a consequence of the following property:
\be
\label{C.10}
\del_t \eta(\vh) + \del_x\big( q_1(\vh) - |\uh|^2 q_2(\vh) \big) \in \{ \text{ compact of }
W^{-1,2}_\mathrm{loc} \},
\ee
where, we recall, $\eta(v)$ is a convex function (the \emph{entropy}), and the \emph{entropy fluxes} $q_{1,2}$ verify $q_1'(v) = \eta'(v) f'(v)$ and $q_2'(v) = \eta'(v) g''(v)$. In practice, one may use the following well-known result to establish \eqref{C.10}: If $1<q<2<r\le\infty$, then 
\be
\label{C.20}
\{\text{ compact of } W^{-1,q}_\loc \} \cap \{ \text{ bounded in } W^{-1,r}_\loc \}
\subset \{ \text{ compact of } W^{-1,2}_\loc \}.
\ee
\begin{lemma}
\label{C-10}
Let $(\uh,\vh)$ be defined by the semidiscrete approximation \eqref{15.50}. Then, the compactness property in \eqref{C.10} is valid. 
\end{lemma}
\begin{proof}
Let $\phi$ be continuous and compactly supported on $\RR \times (0,\infty)$. Let us fix some 
notations. In what follows, $\phi_{j} =\phi(x_j,t)$, $ I_j = [x_j,x_{j+1}]$ and 
$ \overline\phi_{j} = \frac1h \int_{I_j} \phi(y,t) dy$. Let $J = J(h)\in\NN, T>0$ be such that $ \supp\phi
\subset [x_{-J},x_J]\times[0,T]$. Also, we sometimes write, say, $\eta'_{j} := \eta'(\vj)$, and likewise
for the other functions. We view $\vh$ and $\uh$ as piecewise constant functions defined
as $\vh(x) = \vj$ if $x \in I_j$.
We have from \eqref{C.10}
\begin{align}
\label{C.30}
&- \big\langle   \del_t \eta(\vh) + \del_x ( q_1(\vh) - q_2(\vh) |\uh|^2 ), \phi 
\big\rangle_{\mathcal{D}'\times \mathcal{D}}
\\ \notag
&\quad=  -\int_0^\infty h \sum_{j\in\ZZ} \overline\phi_{j} \eta' (\vj) \del_t \vj dt
- \int_0^\infty \sum_{j\in\ZZ} \phi_{j+1} \big( q_1(\vjj) - q_1(\vj) \big) dt 
\\\notag
& \quad\quad+ \int_0^\infty \sum_{j\in\ZZ}\phi_{j+1} \big( q_2(\vjj) |\ujj|^2 - q_2(\vj) |\uj|^2 \big) dt,
\end{align}
where we have just differentiated piecewise constant functions in the sense of distributions.
From the definition of $v^h$ in the scheme \eqref{15.50}, we find, with obvious notation,
\begin{align*}
&-\int_0^\infty h \sum_{j\in\ZZ} \overline\phi_{j} \eta' (\vj) \del_t \vj dt 
= \int_0^\infty  \sum_{j\in\ZZ}\overline\phi_{j}  {\eta'(\vj)}\big( \fjp (\vj, \vjj) + \fjm (\vj, \vjm) \big) dt
\\
&\quad +  \int_0^\infty  \sum_{j\in\ZZ}\overline\phi_{j} \eta'(\vj)\big( \Gjp (\vj, \vjj, |\uj|^2, |\ujj|^2) 
+ \Gjm (\vj, \vjm, |\uj|^2, |\ujm|^2) \big) dt
\\
&\qquad=: A + C.
\end{align*}
Thus, if we define (see \eqref{C.30})
\[
\aligned
B & := - \int_0^\infty \sum_{j\in\ZZ} \phi_{j+1} \big( q_1(\vjj) - q_1(\vj) \big) dt 
\\
D & := \int_0^\infty \sum_{j\in\ZZ}\phi_{j+1} \big( q_2(\vjj) |\ujj|^2 - q_2(\vj) |\uj|^2 \big) \,dt, 
\endaligned
\] 
we get
\be
\label{C.40}
\aligned
&-\big \langle \phi(x,t) \del_t \eta(\vh) + \phi(x,t)\del_x ( q_1(\vh) - q_2(\vh) |\uh|^2 )\big\rangle_{\mathcal{D}'\times \mathcal{D}}
\\
&\quad = A+B+C+D.
\endaligned
\ee
We must therefore estimate $A+B+C+D$ so that the compactness property in \eqref{C.10} is valid.
\subsection*{Estimate of $A+B$}
We begin with
(omitting the arguments of the numerical flux functions whenever possible, so that 
$\fjp = \fjp(\vj, \vjj)$, $\fjm = \fjm(\vj, \vjm)$)
\[
\aligned
A &= \int_0^\infty \sum_{j\in\ZZ} \phi_{j} \etaj' (\fjp + \fjm)\,dt
+ \int_0^\infty \sum_{j\in\ZZ} (\bar\phi_{j} - \phi_{j} ) \etaj' (\fjp + \fjm) \,dt
\\
& =:A_1 + A_2.
\endaligned
\]
First, using the conservation property in \eqref{15.30} and summing by parts, we find
\[
\aligned
A_1 &= \int_0^\infty \sum_{j\in\ZZ}\phi_{j} \etaj'( \fjp(\vj,\vjj) - f_{j-1,+} (\vjm, \vj))\,dt
\\
& = -\int_0^\infty \sum_{j\in\ZZ} (\phi_{j+1} \etajj' - \phi_{j} \etaj' ) \fjp(\vj,\vjj)\,dt
\\
& = -\int_0^\infty \sum_{j\in\ZZ} (\phi_{j+1} - \phi_{j}) \etaj'  \fjp \,dt
-\int_0^\infty \sum_{j\in\ZZ} \phi_{j+1}( \etajj' -  \etaj' ) \fjp\,dt
\endaligned
\]
and so
\[
\aligned
A_1 & = -\int_0^\infty \sum_{j\in\ZZ} (\phi_{j+1} - \phi_{j}) \etaj'  \fjp \,dt 
-\int_0^\infty \sum_{j\in\ZZ} \phi_{j+1}\int_{\vj}^{\vjj} \eta''(v) \fjp \, dv\,dt
\\
& =: A_{11} + A_{12}.
\endaligned
\]
Next, we have
\[
\aligned
B &= -\int_0^\infty \sum_{j\in\ZZ} \phi_{j+1} \int_{\vj}^{\vjj} \eta'(v) f'(v) \,dv\,dt
\\
& = \int_0^\infty \sum_{j\in\ZZ} \phi_{j+1} \int_{\vj}^{\vjj} \eta''(v) f'(v)\, dv\,dt
- \int_0^\infty \sum_{j\in\ZZ} \phi_{j+1} ( \etajj' \fjj - \etaj' \fj ) \,dt 
\\
&=: B_1 + B_2.
\endaligned
\]
We now write
\be
\label{AB}
\aligned
A+B = T_1 + T_2 + T_3, 
\endaligned
\ee
with 
\[
\aligned
T_1 = A_{12} + B_1, \qquad T_2 = A_{11} + B_2 \qquad T_3 = A_2
\endaligned
\]
and we estimate each term in turn. We have immediately
\[
\aligned
|T_1| \le \| \phi \|_\infty \int_0^T \sum_{j\in\ZZ} \int_{\vj}^{\vjj} \eta''(v) ( f(v) - \fjp) \,dv \,dt,
\endaligned
\]
since by Lemma \ref{L-10} the quantity inside the integral is non-negative. Thus, 
by the viscosity estimate \eqref{UB.28} we see that
$T_1$ is bounded in the space of measures (on the support of $\phi$), which is compactly embedded in 
$W^{-1,q}(\supp \phi)$ for $q \in [1,2)$. This proves the compactness property \eqref{C.10} for $T_1$.
Note also that from Lemma \ref{L-10} we have
\be
\label{T1}
T_1 \ge 0 \text{ if }\phi\ge 0,
\ee 
which will be of use later.

We now state a  technical result \cite[Lemma 4.5]{EGH}, which is essential in our
convergence analysis. Adapted to our setting, it states that if $b:\RR \to \RR$ is a monotone function, then 
for some $C>0$, we have the bound
\[
\aligned
\big( b(p) - b(q) \big)^2 \le C \Big| \int_p^q b(s) - b(p) \,ds \Big|.
\endaligned
\]
Since $\fjp(p,p) = f(p)$, this gives
\[
\aligned
\big( f(p) - \fjp(p,q) \big)^2 \le C \Big|\int_p^q \fjp(p,s) - \fjp(p,q) \,ds \Big|.
\endaligned
\]
But the monotonicity of the scheme implies that the integral inside the absolute value is actually
non-negative, and so
\be
\label{EGH}
\aligned
\big( f(p) - \fjp(p,q) \big)^2
&\le C \int_p^q \fjp(p,s) - \fjp(p,q) \,ds 
\\
&\le  C \int_p^q f(s) - \fjp(p,q) \,ds,
\endaligned
\ee
again by monotonicity. Similarly, we find
\be
\label{EGH-1}
\aligned
\big( f(q) - \fjp(p,q) \big)^2
&\le  C \int_p^q f(s) - \fjp(p,q) \,ds.
\endaligned
\ee
With an entirely similar proof, we get also for $\xi = p,q$
\be
\label{EGH-G}
\aligned
\big( - g'(\xi) - \bGjp(p,q) \big)^2
&\le  C \int_p^q -g'(s) - \bGjp(p,q) \,ds.
\endaligned
\ee

With this result in hand, we now treat the term $T_2$, using again the viscosity estimate \eqref{UB.28}. 
First, some straightforward rearranging gives
\be
\label{T2}
\aligned
T_2 &= -\int_0^\infty \sum_{j\in\ZZ} \etaj'  (\phi_{j+1} - \phi_{j})  (\fjp - \fj)\,dt. 
\endaligned
\ee
Suppose now that $\phi$ is 
H\"older continuous with exponent $\alpha \in (1/2, 1)$ and recall that 
$\supp \phi$ is compact, so that in particular the sum in $j$ takes place over a finite set $J$.
We find, using the uniform bound for $\vh$ \eqref{UB.20}, the estimate \eqref{EGH} and the property
$\sum_{j\in J}h = \sum_{j\in J}\int_{I_j} 1 \,dx = \int_{x_{-J}}^{x_J} 1\, dt \le C|\supp \phi(t,\cdot)|$,
\[
\aligned
| T_2| &\le C\|\eta'(\vh)\|_\infty h^\alpha \|\phi\|_{0,\alpha} \int_0^T \sum_{j\in J}
| \fjp - \fj| \,dt
\\
&\le C h^\alpha \|\phi\|_{0,\alpha} \int_0^T \Big(\sum_{j\in J} h \Big)^{1/2}
\Big(\sum_{j\in J} (\fjp - \fj)^2 \Big)^{1/2} h^{-1/2} \,dt
\\
&\le C h^{\alpha-1/2} \|\phi\|_{0,\alpha} \int_0^T \Big(\sum_{j\in J} \int_{\vj}^{\vjj}
 f(v) - \fjp \, dv \Big)^{1/2} dt.
\endaligned
\]
Therefore we may use the viscosity estimate \eqref{UB.28} to get
\be
\label{T2-1}
\aligned
|T_2| \le C h^{\alpha-1/2} \|\phi\|_{0,\alpha} a(T), 
\endaligned
\ee
where $a(t)$ is a bounded function. Since $W^{1,q'} \subset C^{0,\alpha},$ with compact embedding, for $q' \ge 2/(1-\alpha)>4$ (and thus $q\in(1,4/3)$), we see that this term is compact in $W^{-1,q}_\loc$ for $q \in (1,4/3) \subset (1,2)$.

The term $T_3$ is treated in the same way, by putting
\[ \fjp(\vj, \vjj) + \fjm(\vj, \vjm) = \fjp(\vj,\vjj) - f(\vj) - f_{j-1,+}(\vjm, \vj) + f(\vj),\]
using \eqref{EGH-1}, and observing that $|\bar\phi_j - \phi_j|\le h^\alpha \|\phi\|_{0,\alpha}$. Thus,
the terms $A+B$ in \eqref{C.40} have the desired compactness property.
\subsection*{The terms $C+D$}
We now turn to the remaining terms in \eqref{C.40}. These are non-homogenous terms involving
the solution $\uh$ to the Sch\-r\"o\-din\-ger equation, and thus present additional difficulties.
Recall that $C$ is defined after \eqref{C.30}. First, we must put $\phi_j$ in place of $\bar\phi_{j}$.
We have
\[
\aligned
C &= \int_0^\infty \sum_{j\in\ZZ} \phi_{j} \etaj' ( \Gjp + \Gjm) \,dt+ \int_0^\infty \sum_{j\in\ZZ}
(\overline\phi_{j} - \phi_{j}) \etaj' ( \Gjp + \Gjm) \,dt
\\ 
& =: C_1 + C_2.
\endaligned
\]
Now, from the conservation property of the numerical flux in \eqref{15.40},
\[
\aligned
C_1 &= \int_0^\infty \sum_{j\in\ZZ} \phi_{j} \etaj'  \big( \Gjp(\vj,\vjj,|\uj|^2, |\ujj|^2) 
- G_{j-1,+}(\vjm,\vj,|\ujm|^2, |\uj|^2) \big) \,dt
\\
& = - \int_0^\infty \sum_{j\in\ZZ} (\phi_{j+1} \etajj' - \phi_{j} \etaj') \Gjp(\vj,\vjj,|\uj|^2, |\ujj|^2) \,dt
\\
& = -  \int_0^\infty \sum_{j\in\ZZ} (\phi_{j+1}  - \phi_{j}) \etaj' \Gjp \,dt
-  \int_0^\infty \sum_{j\in\ZZ} \phi_{j+1}( \etajj' -  \etaj') \Gjp \,dt
\endaligned
\]
and so
\[
\aligned
C_1  & = -  \int_0^\infty \sum_{j\in\ZZ} (\phi_{j+1}  - \phi_{j}) \etaj' \Gjp \,dt
-  \int_0^\infty \sum_{j\in\ZZ} \phi_{j+1} \int_{\vj}^{\vjj} \eta''(v) \Gjp \,dv\,dt
\\
& =: C_{11} + C_{12}.
\endaligned
\]
On the other hand,
\[
\aligned
D &= \int_0^\infty \sum_{j\in\ZZ}\phi_{j+1} \big( q_2(\vjj) |\ujj|^2 - q_2(\vj) |\uj|^2 \big) \,dt
\\
& = \int_0^\infty \sum_{j\in\ZZ}\phi_{j+1} |\ujj|^2 \big( q_2(\vjj)  - q_2(\vj) \big) \,dt
\\
&\qquad+ \int_0^\infty \sum_{j\in\ZZ}\phi_{j+1}  q_2(\vj)\big( |\ujj|^2 - |\uj|^2 \big) \,dt.
\endaligned
\]
But since $q_2' = \eta' g''$, we find
\[
\aligned
&\sum_{j\in\ZZ}\phi_{j+1} |\ujj|^2 \big( q_2(\vjj)  - q_2(\vj) \big) 
\\
&\qquad =  \sum_{j\in\ZZ}\phi_{j+1} |\ujj|^2 \int_{\vj}^{\vjj} \eta'(v) g''(v) \,dv 
\\
&= -  \sum_{j\in\ZZ}\phi_{j+1} |\ujj|^2 \int_{\vj}^{\vjj} \eta''(v) g'(v) \,dv 
+  \sum_{j\in\ZZ} \phi_{j+1} |\ujj|^2 \big( \etajj' \gjj  - \etaj' \gj \big) 
\endaligned
\]
and so
\[
\aligned
D & = - \int_0^\infty \sum_{j\in\ZZ}\phi_{j+1} |\ujj|^2 \int_{\vj}^{\vjj} \eta''(v) g'(v) \,dv \,dt 
\\
&\qquad + \int_0^\infty \sum_{j\in\ZZ} \phi_{j+1} |\ujj|^2 \big( \etajj' \gjj  - \etaj' \gj \big) \,dt
\\
&\qquad \qquad +\int_0^\infty \sum_{j\in\ZZ}\phi_{j+1}  q_2(\vj)\big( |\ujj|^2 - |\uj|^2 \big) \,dt
\\
&=: D_1 + D_2 + D_3.
\endaligned
\]
We now write
\be
\label{CD}
\aligned
C + D = S_1 + S_2 + D_3 + C_2,
\endaligned
\ee
with
\[
\aligned
S_1 := D_1 + C_{12}, \qquad S_2 := D_2 + C_{11}
\endaligned
\]
and still (adding and subtracting)
\be
\label{S1}
\aligned
S_1 &= -\int_0^\infty \sum_{j\in\ZZ}\phi_{j+1} |\ujj|^2 \int_{\vj}^{\vjj} \eta''(v) 
\big( g'(v) + \bGjp(\vj,\vjj) \big) \,dv \,dt 
\\
& + \int_0^\infty \sum_{j\in\ZZ}\phi_{j+1}  \int_{\vj}^{\vjj} \eta''(v) 
\big( |\ujj|^2 \bGjp(\vj, \vjj) - \Gjp \big) \,dv \,dt 
\\
&=: S_{11} + S_{12}.
\endaligned
\ee
Thus we have succeeded in bringing out the viscosity of the scheme in $S_{11}$.
We now bound these two terms. We have
\[
\aligned
|S_{11}| & \le \| \phi\|_\infty \int_0^T \sum_{j\in\ZZ} |\ujj|^2 \int_{\vj}^{\vjj} \eta''(v) 
\big(- g'(v) - \bGjp(\vj,\vjj) \big) \,dv \,dt 
\\
&\le C \|\phi\|_\infty, 
\endaligned
\]
in view of the viscosity bound \eqref{UB.28} (the bound \eqref{UB.28} actually involves
$|\uj|^2$ and not $|\ujj|^2$, but an examination of the proof quickly gives the same result with
$|\ujj|^2$). Thus, $S_{11}$ is bounded in the space of
measures, which as we have seen in the treatment of the term $T_1$ above, is sufficient 
for our purposes. Additionally, we find from Lemma~\ref{L-10}
\be
\label{S11}
S_{11} \ge 0 \text{ if } \phi \ge 0,
\ee
which will be useful later.

Next, we have
\[
\aligned
S_{12} &\le C\|\phi\|_\infty \int_0^T \sum_{j\in\ZZ} \big| |\ujj|^2 - |\uj|^2 \big| \,dt
\le C \|\phi\|_\infty \int_0^T \|D_+ \uh\|_2 dt,
\endaligned
\]
by using the uniform $L^\infty$ bound on $\vh$, the bound \eqref{UB.101} and the smoothness of $\Gjp$ (recall that 
$\Gjp(a,b,u,u) = u\bGjp(a,b)$). In view of the bound on $\|D_+\uh\|_2$, \eqref{UB.27}, this term is also 
bounded in the space of measures.

Consider now the term $S_2 = D_2 + C_{11}$ in \eqref{CD}. Some easy rearranging, addition and 
subtraction give
\be
\label{S2}
\aligned
S_2 &= \int_0^\infty \sum_{j\in\ZZ} - (\phi_{j+1} - \phi_{j}) |\uj|^2 \etaj' \big( \gj + \bGjp \big) 
\\
& \qquad 
+ (\phi_{j+1} - \phi_{j})  \etaj' \big( |\uj|^2 \bGjp - \Gjp \big) 
- \phi_{j+1} \etaj' \gj (|\ujj|^2 - |\uj|^2) \,dt.
\endaligned
\ee
Now the first of these three terms is treated in the same way as the term \eqref{T2} above,
using the viscosity estimate \eqref{UB.28}, the Sobolev inequality
\eqref{M.10}, and the technical estimate \eqref{EGH-G}, giving
\[
\aligned
&\int_0^\infty \sum_{j\in\ZZ} - (\phi_{j+1} - \phi_{j}) |\uj|^2 \etaj' \big( \gj + \bGjp \big) 
\\
&\qquad \le  h^\alpha \|\phi\|_{0,\alpha} \|\eta'(\vh)\|_\infty  \int_0^T\|\uh(t)\|_{\infty} \sum_{j\in J}
|\uj| |\gj + \bGjp | \,dt
\\
&\qquad \le C h^\alpha \|\phi\|_{0,\alpha}  \int_0^T \sqrt{a(t)} \Big(\sum_{j\in J} h \Big)^{1/2}
\Big(\sum_{j\in J} |\uj|^2 (\gj + \bGjp )^2 \Big)^{1/2} h^{-1/2} \,dt
\\
&\qquad \le C h^{\alpha-1/2} \|\phi\|_{0,\alpha} \int_0^T\sqrt{a(t)} \Big(\sum_{j\in J}|\uj|^2 \int_{\vj}^{\vjj} 
-g'(v) - \bGjp \, dv \Big)^{1/2} dt
\\
&\qquad \le C h^{\alpha-1/2} \|\phi\|_{0,\alpha} A(T)
\endaligned
\]
for some locally bounded function $A(t)$. This gives compactness 
in $W^{-1,q}_\loc$ for $q \in (1,4/3) \subset (1,2)$. 
The second and third terms may be treated exactly as $S_{12}$ above, and are thus bounded
in the space of measures and so have the desired compactness. Note also that the term $D_3$ has
precisely the same form and so yields to the same analysis.

It only remains to estimate the term
\be
\label{C2}
\aligned
C_2 =\int_0^\infty \sum_{j\in\ZZ}
(\overline\phi_{j} - \phi_{j}) \etaj' ( \Gjp + \Gjm) \,dt
\endaligned
\ee
which we only sketch, since no new techniques are necessary: just add and subtract appropriately
to obtain terms with $(\overline\phi_{j} - \phi_{j})|\uj|^2(\gj + \bGjp)$ and 
$(\overline\phi_{j} - \phi_{j}) (\Gjp - |\uj|^2 \bGjp )$, which fall into the
cases treated before.
This completes the proof of Lemma~\ref{C-10}.\end{proof}

\begin{lemma}
\label{conv-1}
Let $(\uh, \vh)$ be defined by the scheme \eqref{15.50}. Then,
there exist functions 
\[
\aligned
u \in L^\infty_\loc( [0,\infty) ; H^1(\RR)), \qquad v \in L^\infty(\RR \times [0, \infty))
\endaligned
\]
such that (for a subsequence at least)
\[
\aligned
&\vh \to v  \quad \uh \to u \quad\text{in} \quad L^1_\loc( \RR \times [0, \infty) )
\endaligned
\]
and 
\[
D_+ \uh \ws \del_x u \quad\text{in} \quad L^\infty([0,\infty) ; L^2(\RR)).
\]
\end{lemma}
\begin{proof}
We have proved the compactness property \eqref{C.10}. Thus, according to the compensated 
compactness method \cite{Tartar}, there is a function $v \in L^1_\loc( \RR \times [0, \infty) )$ (and 
thus in $L^\infty( \RR \times [0, \infty) ),$ by the uniform bound \eqref{UB.20})
such that $\vh \to v$ in $ L^1_\loc( \RR \times [0, \infty) )$.

Let us now analyze the compactness of $\uh$. For this, it will be useful to define the
piecewise interpolators $P_1$ and $P_0$. The piecewise linear interpolator
$P_1 \uh$ is the unique continuous, 
piecewise linear function (on each interval $I_j = [x_j, x_{j+1})$) such that $P_1 \uh (x_j) = \uj$, 
and the piecewise constant interpolator is defined as $P_0 \uh(x) = \uj $ for $x \in I_j$.
In view of the estimate $\| D_+\uh\|_{2} \le a(t)$ and the fact that
\[
\| D_+ \uh \|_{2} \equiv \| P_0 D_+ \uh \|_2 = \| \del_x P_1\uh \|_2,
\]
we see that there exists a function $u \in L^\infty_\loc( [0, \infty); H^1(\RR))$
such that $P_1 \uh \to u$ 
in $L^1_\loc([0,\infty) \times\RR)$  and
$\del_x P_1 \uh \ws \del_x u $ in $L^\infty([0,\infty);L^2(\RR))$. Note that this proves the last
assertion of the Lemma, since $D_+ \uh = \del_x P_1 \uh$. With this in mind, we wish
to prove that $\uh \equiv P_0 \uh \to u$ in $L^1_\loc(\RR \times[0,\infty))$; 
For this, let $\Omega$ be a bounded set of $\RR \times[0,\infty)$. We have
\[
\aligned
\| P_0 \uh - u \|_{L^1(\Omega)} &\le
\| P_0 \uh -  P_1 \uh \|_{L^1(\Omega)} 
+ \| P_1 \uh - u \|_{L^1(\Omega)}.
\endaligned
\]
We have just seen that the second term above vanishes as $h\to 0$. For the first term, 
note that since $\Omega$ is compact, there are some $M, T \in \RR$ independent of $h$ such that
\[
\aligned
\| P_0 \uh -  P_1 \uh \|_{L^1(\Omega)}  \le \int_0^T \sum_{|j|< M/h} \int_{x_j}^{x_{j+1}}
|P_0 \uh(x,t) -  P_1 \uh (x,t)| \, dx dt.
\endaligned
\]
Since
\[
\aligned
\int_{x_j}^{x_{j+1}} |P_1\uh - P_0 \uh| dx = \frac{h}2 |\ujj - \uj|,
\endaligned
\]
we obtain from the estimate $\| D_+\uh\|_{2} \le a(t)$ in Proposition~\ref{UB-20} and
Young's inequality that
\[
\aligned
\| P_0 \uh -  P_1 \uh \|_{L^1(\Omega)}  &\le 
\int_0^T \sum_{|j|< M/h} \frac{h}2 |\ujj - \uj| \, dt
\\
&\le C \int_0^T \Big( \sum_{|j|< M/h} h^2 + \sum_{j\in\ZZ} |\ujj - \uj|^2 \Big) \, dt
\\
&\le C \int_0^T h (M + a(t)) \,dt
\endaligned
\]
which tends to zero with $h$. This proves that $P_0 \uh \to u $ in $L^1_\loc(\RR \times[0,\infty))$
and completes the proof of Lemma~\ref{conv-1}.
\end{proof}

\proof[Proof of Theorem \ref{20-10}]
It only remains to check that the (strong) limit $(u,v)$ of $(\uh, \vh)$ given by 
Lemma \ref{conv-1}
is the unique entropy solution to the Cauchy problem \eqref{I.10}--\eqref{I.25}. 

First, consider the 
Sch\-r\"o\-din\-ger equation \eqref{I.10}. Taking the discrete equation \eqref{15.51}, multiplying
by a test function $\theta$ as in Definition \ref{S-05} and integrating gives
\[
\aligned
\iint_{\RR\times[0,\infty)} -i\uh\del_t \theta  + \theta D_h^2 \uh - (|\uh|^2 \uh \theta
+ g(\vh) \uh )\theta\, dx dt = 0.
\endaligned
\]
From the convergence properties in Lemma \ref{conv-1}, we see that all terms except 
the one with $D_h^2  \uh$ converge to the corresponding ones in the equation \eqref{I.10}.
Regarding that term, we easily find from the definition of the discrete derivatives
\[
\aligned
\iint_{\RR\times[0,\infty)} \theta D_h^2 \uh \,dx dt
= - \iint_{\RR\times[0,\infty)} D_+ \uh \theta_h' \, dxdt,
\endaligned
\]
where $\theta_h'(x,t) := ( \medint_{I_{j+1}} \theta(x,t) \,  dx - \medint_{I_{j}} \theta(x,t) \, dx)/h.$
Since $\theta_h' $ converges strongly to $\del_x \theta$, Lemma~\ref{conv-1} ensures that
this term also converges to the corresponding one in \eqref{I.10}.

Next, we wish to prove that $v$ is a solution to the conservation law \eqref{I.20}, in the
sense of Definition~\ref{S-05}. To this end,
we must prove that for every non-negative smooth function $\phi$ with compact support in 
$ \RR \times [0,T) $, the inequality
\be
\label{C.45}
\aligned
& \iint_{\RR\times[0,\infty)}  \eta(v) \del_t \phi  +  ( q_1(v) - q_2(v) |u|^2 ) \del_x \phi
\\ 
&\quad+  \big( \eta'(v) g'(v) - q_2(v) \big) \del_x (|u|^2) \phi  \, dx dt + \int_\RR \eta(v_0(x)) \phi(x,0) \,dx \ge 0
\endaligned
\ee
holds. What we have is (cf. \eqref{C.30})
\begin{align}
\label{C.50}
&\iint_{\RR\times[0,\infty)}    \eta(\vh) \del_t\phi + ( q_1(\vh) - q_2(\vh) |\uh|^2 ) 
\del_x \phi\, dx dt
\\ \notag
&=  -\int_0^\infty h \sum_{j\in\ZZ} \overline\phi_{j} \eta' (\vj) \del_t \vj dt
- \int_0^\infty \sum_{j\in\ZZ} \phi_{j+1} \big( q_1(\vjj) - q_1(\vj) \big) dt 
\\\notag
& \quad+ \int_0^\infty \sum_{j\in\ZZ}\phi_{j+1} \big( q_2(\vjj) |\ujj|^2 - q_2(\vj) |\uj|^2 \big) dt 
- \int_\RR \eta(v^h_0(x)) \phi(x,0) \,dx,
\end{align}
and so
\be
\label{C.60}
\aligned
&\iint_{\RR\times[0,\infty)}  \eta(\vh)\del_t\phi + ( q_1(\vh) - q_2(\vh) |\uh|^2 ) \del_x  \phi\, dx dt
\\
& = A+B+C+D -  \int_\RR \eta(v^h_0(x)) \phi(x,0) \,dx,
\endaligned
\ee
where $A,B,C$ and $D$ are defined before \eqref{C.40}.

From the strong convergence of $(\uh,\vh)$ to $(u,v)$ given in Lemma~\ref{conv-1},
we see that the left-hand side of \eqref{C.60} converges to the
first line of \eqref{C.45}. Since (by assumption) $v^h_0 \to v_0$ strongly,
it it only necessary to check that 
\be
\label{C.70}
\aligned
\lim_{h\to0} (A+B+C+D) \ge 
-\iint_{\RR\times[0,\infty)} 
\phi(x,t) \big( \eta'(v) g'(v) - q_2(v) \big) \del_x |u|^2   \, dx dt
\endaligned
\ee
for \eqref{C.45} to hold.

Consider the terms $A+B$. From the decomposition \eqref{AB}, the positivity 
property \eqref{T1}, the estimate \eqref{T2-1} and the discussion after \eqref{T2-1},
we see that $\lim_{h\to 0} A+B \ge 0 .$

Let us now analyze the terms contained in $C+D$, see \eqref{CD}.
First, we claim that $\lim_{h\to 0} S_1 \ge 0.$ Recalling \eqref{S1}, the property \eqref{S11} shows that
$S_{11}$ has a sign. So
we must prove that 
\[
\aligned
S_{12} \equiv \int_0^\infty \sum_{j\in\ZZ}\phi_{j+1}  \int_{\vj}^{\vjj} \eta''(v) 
\big( |\ujj|^2 \bGjp(\vj, \vjj) - \Gjp \big) \,dv \,dt \to 0
\endaligned
\]
as $h\to0$. We have
\[
\aligned
S_{12} &\le C \|\phi\|_\infty \int_0^T \sum_{|j|\le J} \big| \etajj'  - \etaj' \big|
\big| |\ujj|^2 - |\uj|^2 \big| \,dt
\\
& \le C \Big( \int_0^T \sum_{|j|\le J}h \big| \vjj  - \vj \big|^2 \,dt \Big)^{1/2}
\Big( \int_0^T \sum_{|j|\le J}\frac1h \big| |\ujj|^2 - |\uj|^2 \big|^2 \,dt \Big)^{1/2}
\endaligned
\]
and, from \eqref{UB.27},\eqref{M.10},
\[
\aligned
\sum_{|j|\le J}\frac1h \big| |\ujj|^2 - |\uj|^2 \big|^2  & \le \|\uh\|_\infty^2
\sum_{j\in\ZZ} \frac1h |\ujj - \uj|^2 \le C(T).
\endaligned
\]
Therefore, $S_{12} \to 0$ if 
\[
\int_0^T \sum_{|j|\le J}h \big| \vjj  - \vj \big|^2 \,dt \to 0
\]
or, in particular, if
\be
\label{S12}
\iint_{\Omega} \big| \vh_+  - \vh \big|^2 \,dx\,dt \equiv \| \vh - \vh_+\|_{L^2(\Omega)}^2 \to 0,
\ee
where we have set $v_+(x) = v(x+h)$ and $\Omega = (0,T) \times \RR$.
A standard result in compactness \cite[p.~4]{Evans} implies that if $\| \vh\|_{L^2(\Omega)} 
= \| \vh_+\|_{L^2(\Omega)}$, and if 
$\vh \to v$ in $L^2(\Omega)$, then $\vh_+ \rightharpoonup v$ actually gives
$\vh_+ \to v$, and so \eqref{S12} will hold. Therefore, we only have to prove the weak convergence
$\vh_+ \rightharpoonup v$ in $L^2(\Omega)$. But if $\varphi \in L^2(\Omega)$ then
\[
\aligned
\iint_\Omega & \varphi (t,x) (\vh(t,x+h) - v(t,x) )\,dx \,dt 
\\& =
\iint_\Omega \varphi (t,x) (\vh(t,x+h) - \vh(t,x) ) +
 \varphi (t,x) (\vh(t,x) - v(t,x) )\,dx \,dt
\\
& = \iint_\Omega (\varphi (t,x - h) - \varphi(x)) \vh(t,x) 
+  \varphi (t,x) (\vh(t,x) - v(t,x) )\,dx \,dt
\\
& \le \|\varphi_+ - \varphi \|_2 \| \vh_+\|_2 + \| \varphi\|_2 \| \vh - v\|_2 \to 0.
\endaligned
\]
This proves \eqref{S12} and so the term $S_{12}$ goes to zero with $h$.

Next, we estimate the term $S_2$ given in \eqref{S2}. We find
\[
\aligned
S_2 &= \int_0^\infty \sum_{j\in\ZZ} - (\phi_{j+1} - \phi_{j}) |\uj|^2 \etaj' \big( \gj + \bGjp \big) 
\\
& \qquad 
+ (\phi_{j+1} - \phi_{j})  \etaj' \big( |\uj|^2 \bGjp - \Gjp \big) 
- \phi_{j+1} \etaj' \gj (|\ujj|^2 - |\uj|^2) \,dt
\\
& =: S_{21} + S_{22} + S_{23}.
\endaligned
\]
Thus
\[
\aligned
S_{21} \le h \|\phi' \|_\infty \int_0^T \sum_{|j| \le J} |\uj|^2 \etaj' \big(- \gj - \bGjp \big) ,
\endaligned
\]
and so the same computation as the one after \eqref{S2} applies, giving 
\[
S_{21} \to 0.
\]
Next, using the estimates \eqref{UB.101} and \eqref{UB.27},
\[
\aligned
S_{22} &\le C h \|\phi'\|_\infty \int_0^T \sum_{j\in\ZZ} \big| |\ujj|^2 - |\uj|^2 \big| \,dt
\\
&\le C(T) h \|\phi'\|_\infty \to 0
\endaligned
\]
as $h\to 0$. Finally, from Lemma~\ref{conv-1} and from
$D_+ |\uh|^2 \rightharpoonup \del_x |u|^2$ (which is a consequence of the same lemma),
we see that 
\[
\aligned
S_{23} \to - \iint_\Omega \phi(x,t) \eta'(v) g'(v) \del_x|u|^2 \,dx \,dt .
\endaligned
\]
It only remains to analyze the terms $D_3$ and $C_2$ in \eqref{CD}. We have,
again from Lemma~\ref{conv-1},
\[
\aligned
D_3 &= \int_0^\infty \sum_{j\in\ZZ}\phi_{j+1}  q_2(\vj)\big( |\ujj|^2 - |\uj|^2 \big) \,dt
\\
& \to \iint_\Omega \phi(x,t) q_2(v) \del_x |u|^2 \,dx \,dt
\endaligned
\]
as $h\to 0$. For the term $C_2$, we easily reduce it to the cases already treated by proceeding as in the
comments after \eqref{C2}, giving $C_2 \to 0$.
This completes the proof of Theorem~\ref{20-10}. \endproof


\section{Numerical experiments}
\label{NE}
In this section, we present some numerical computations illustrating our results.
\subsection{A fully discrete scheme}
To implement numerically the scheme \eqref{15.50}, we choose a semi-implicit Crank--Nicholson scheme
for the Sch\-r\"o\-din\-ger equation \eqref{15.51}, with a Newton iteration for the nonlinear term
$|u|^2u$, coupled to a semi-implicit Lax--Friedrichs scheme (see \eqref{LF-01}) for the conservation law \eqref{15.52}.
More precisely, given a spatial mesh size $h$ and a time step $\tau$,
we consider the following algorithm:
\be
\label{scheme}
\aligned
& i \frac1\tau( \unnj - \unj) +  \frac1{2h^2}\big( \unnjj + \unjj - 2(\unnj + \unj) 
+ \unnjm + \unjm \big)
\\
&\qquad = \big| \frac12(\unnj + \unj) \big|^2 \frac12(\unnj + \unj) + g(\vnj) \frac12(\unnj + \unj),
\\
& \frac1\tau (\vnnj - \vnj) = - \frac1{2h}(f(\vnjj) - f(\vnjm)) 
\\
&\qquad
+ \frac1{2h}
(g'(\vnjj) |\unjj|^2 - g'(\vnjm)|\unjm|^2)
 + \frac1{2\lambda h}(\vnnjj -2\vnnj + \vnnjm)
\\
&\qquad + \frac1{2\gamma h}
\big[ (\vnnjj - \vnnj) \frac12( |\unj|^2 + |\unjj|^2) - (\vnnj - \vnnjm) \frac12( |\unj|^2 + |\unjm|^2)
\big].
\endaligned
\ee
Of course this scheme, being at most first order, serves only to illustrate our results. Indeed, it is
well known that the Lax--Friedrichs scheme suffers from a degradation of solutions due to
numerical viscosity. More accurate 
discretizations would be interesting to implement, for instance, by considering more 
powerful time-stepping methods, or by implicitly coupling the two equations.

As is usual in simulations involving the Sch\-r\"o\-din\-ger equation on the whole line, 
it is reasonable to consider the problem in a sufficiently large domain, and taking care that the
value of the computed solution remains very small (in our case, this value is
usually of the order $10^{-7}$ after a few
thousand iterations) near the boundary.

In all simulation presented below, the computed $L^2$ norm of the solution $\uh$ is 
exactly conserved (up to a very small error committed in its computation), in accordance with the 
$L^2$ norm conservation property \eqref{UB.10}.

\subsection{A test case with linear flux}
As a first example, and to test the accuracy of the numerical scheme \eqref{scheme},
we consider the following problem in which the conservation law has a linear flux,
\be
\label{Laurencot}
\left\{
\aligned
&i\del_t u + \del_{xx} u = k \big( v u + q |u|^2u \big)
\\
&\del_t v + \gamma \del_x v = \delta \del_x(|u|^2)
\endaligned\right.
\ee
whose exact traveling wave solution can be found (see, for instance, \cite{Laurencot}): 
\be
\label{Exact-01}
\aligned
&u(x,t) = e^{i\lambda t} e^{ic(x-ct)/2} \sqrt{\frac{2E}{|\beta|}} \mathop{\mathrm{sech}}\big(
\sqrt{E}(x-ct)\big),
\\
&v(x,t) = a \frac{2E}{|\beta|} {\mathop{\mathrm{sech}}}^2\big(
\sqrt{E}(x-ct)\big),
\endaligned
\ee
where $c,\lambda > 0$ are to be chosen, $E := \lambda -c^2/4$, $a:= \delta/(\gamma-c)$ and 
$\beta = k(a+q)$. We must have $E>0$, $\beta<0$ to ensure
existence of a solution (see \cite{Laurencot} for the details). 
Here, we choose $q = k = \lambda = \gamma =1$, 
$a = -2$ and $c = 3/2$.

Note that here and in the following examples, we have omitted the function $g(v)$ from the 
formulation of the problem. This is because we are taking $g'(v)$ to be the
characteristic function of some large interval $[-M,M]$. In all our numerical experiments,
the value of $|v|$ never gets close to $M$, and thus among the tested initial data,
the original problem can be stated without reference to $g$.
\begin{figure}
\includegraphics[width=\linewidth,keepaspectratio=false]{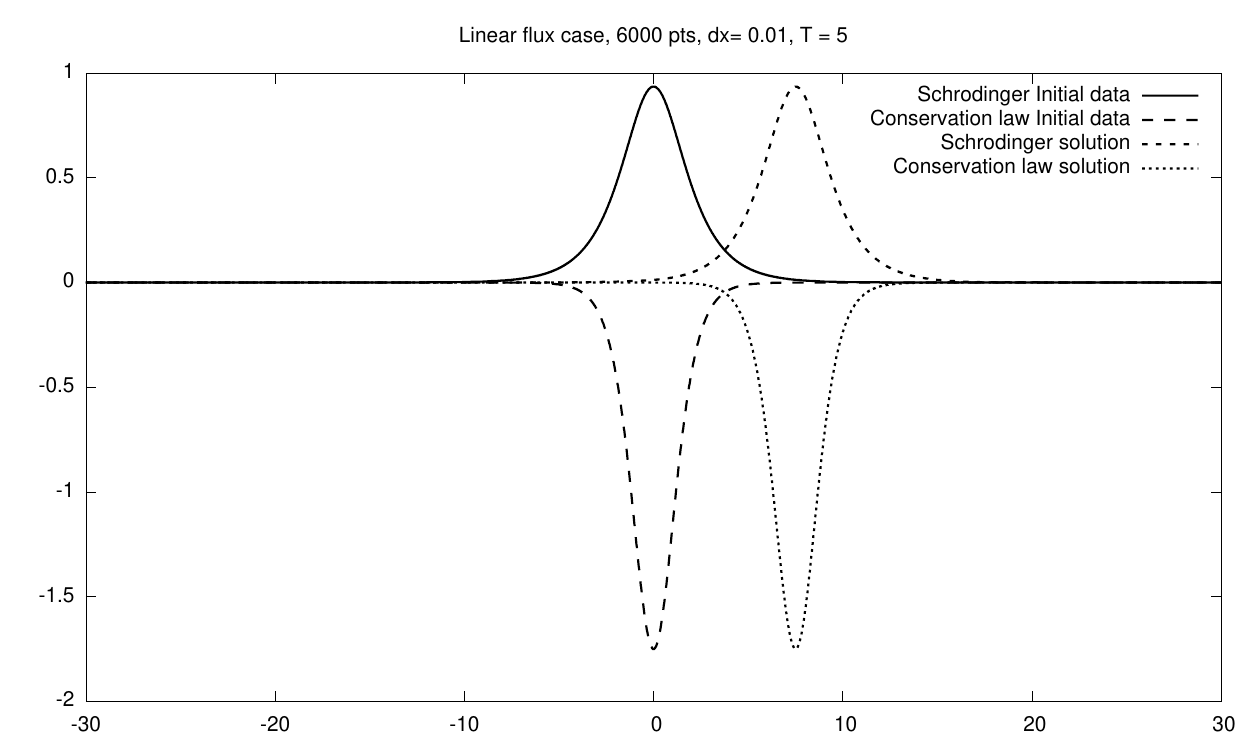}
\caption{Initial data and computed solutions of system \eqref{Laurencot}.}
\label{Linear-01}
\end{figure}
\begin{figure}
\includegraphics[width=\linewidth,keepaspectratio=false]{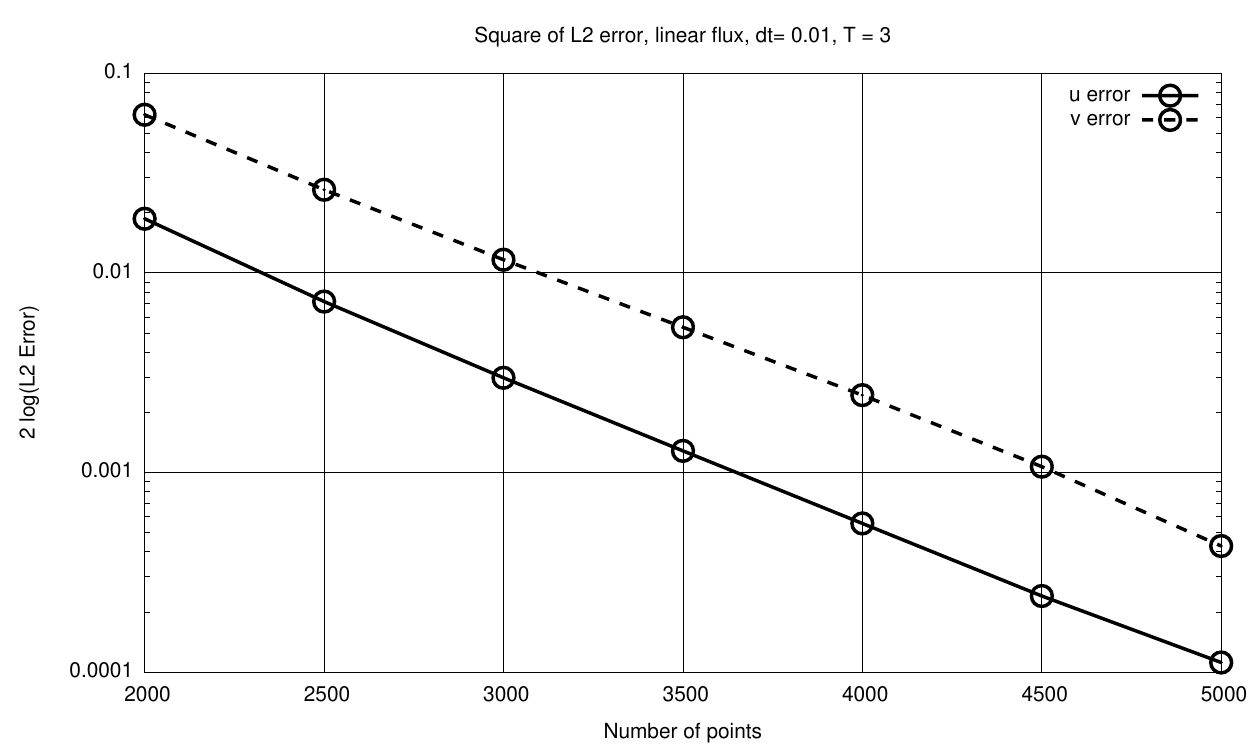}
\caption{$L^2$ error between the 
exact solution \eqref{Exact-01} and the
computed solution, linear flux.}
\label{Linear-02}
\end{figure}

In Figure \ref{Linear-01} we present the initial data $|u_0|$ and $v_0$ corresponding to
\eqref{Exact-01} with $t=0$, 
together with the computed solution $|\uh|$, $\vh$, which can 
be observed to be a traveling wave. The error in the $L^2$ norm between the 
exact solution \eqref{Exact-01} and the
computed solution for various spatial discretizations is presented in Figure~\ref{Linear-02}.

\subsection{A test case with nonlinear flux}
\begin{figure}
\includegraphics[width=\linewidth,keepaspectratio=false]{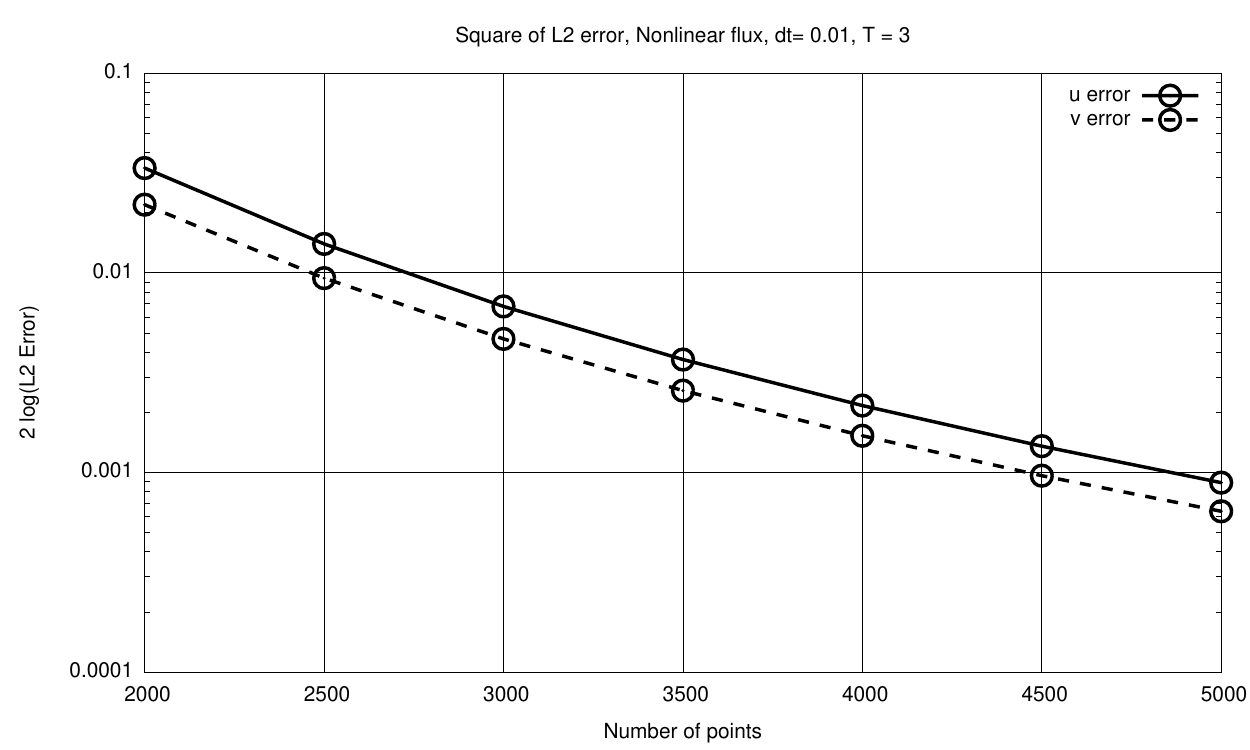}
\caption{$L^2$ error between the 
exact solution and the
computed solution, nonlinear flux.}
\label{Nonlinear-01}
\end{figure}

We now test the scheme \eqref{scheme} on an explicit solution to the problem with 
nonlinear flux $f(v) = v^2$ but in which the first equation is linear in $u$
(we thank J.P. Dias for bringing this solution to our attention):
\[
\left\{
\aligned
&i\del_t u + \del_{xx} u = v u
\\
&\del_t v + \del_x v^2 = \del_x(|u|^2).
\endaligned\right.
\]
The explicit solutions can be found by the following \emph{Ansatz,}
\[
\aligned
&u(x,t) = e^{i b t} r(x),
\\
&v(x,t) = \beta(x).
\endaligned
\]
From the equations we get
\[
\aligned
\beta = -|r| \text{ and } -br + r'' = -r^2,
\endaligned
\]
of which a solution is known to be $r(x) = b(3/2) \mathop{\mathrm{sech}}^2(\sqrt{b}x/2)$. Thus
\[ (u,v) = ( e^{ibt} r(x), -r(x)). \]
In Figure \ref{Nonlinear-01} we present the error in the $L^2$ norm between this exact solution and the
computed solution for various spatial discretizations.

\subsection{The general case}
Here we present some numerical solutions to the full system \eqref{I.10},\eqref{I.20}, 
with $f(v) = 3v^2$, showing 
(the modulus of) a solution of the Sch\-r\"o\-din\-ger equation interacting
with a solution of the conservation law. 

\begin{figure}
\includegraphics[width=\linewidth,keepaspectratio=false]{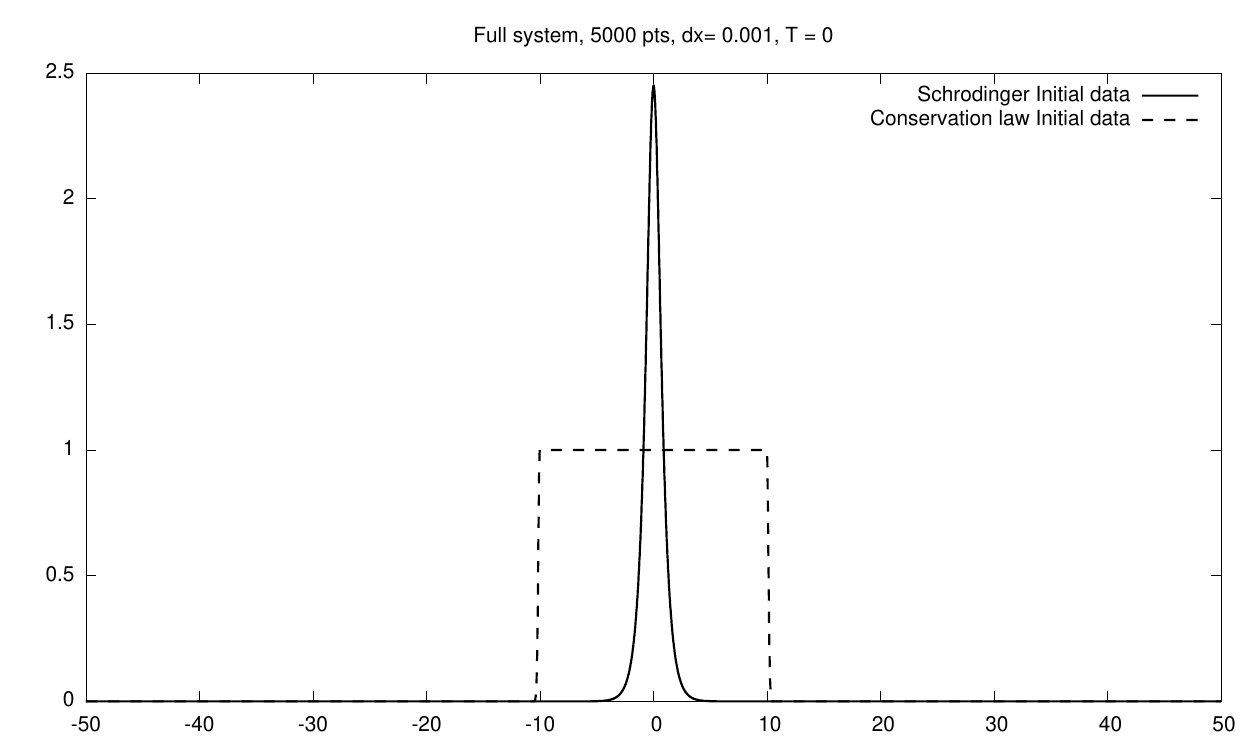}
\includegraphics[width=\linewidth,keepaspectratio=false]{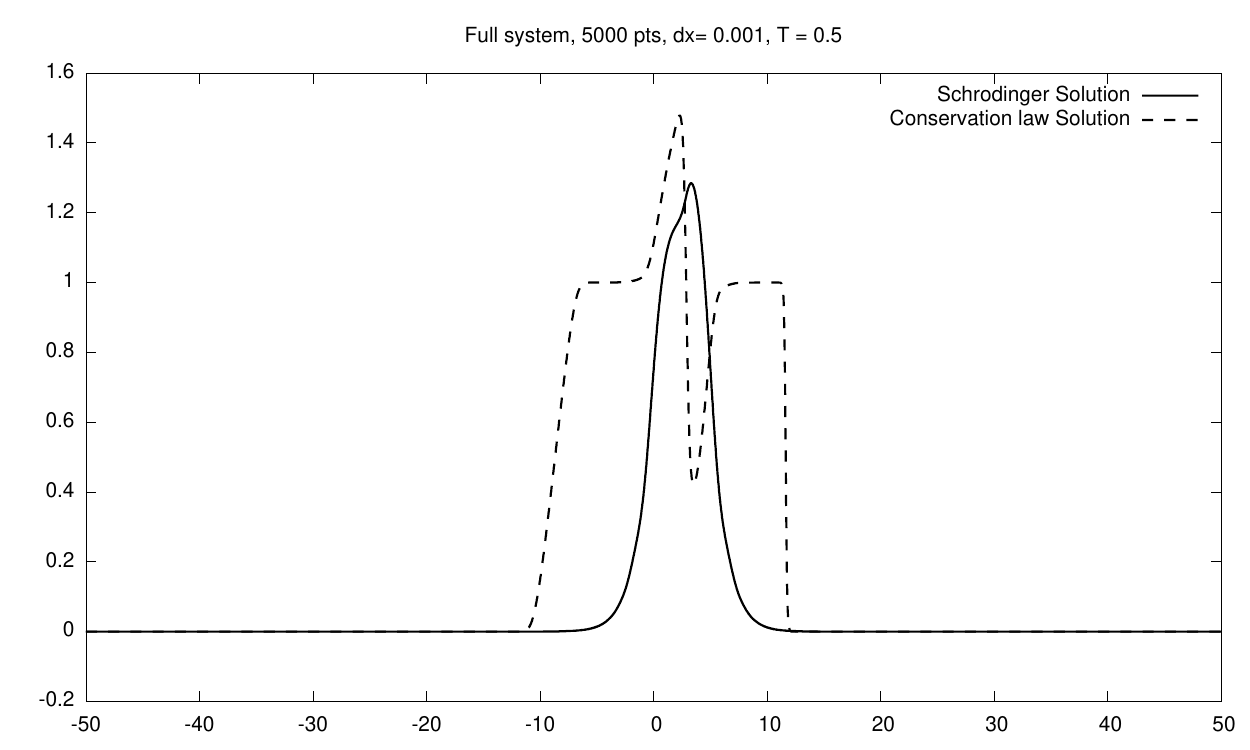}
\caption{Initial data and computed solutions of system \eqref{I.10},\eqref{I.20}.}
\label{FN-01}
\end{figure}

\begin{figure}
\includegraphics[width=\linewidth,keepaspectratio=false]{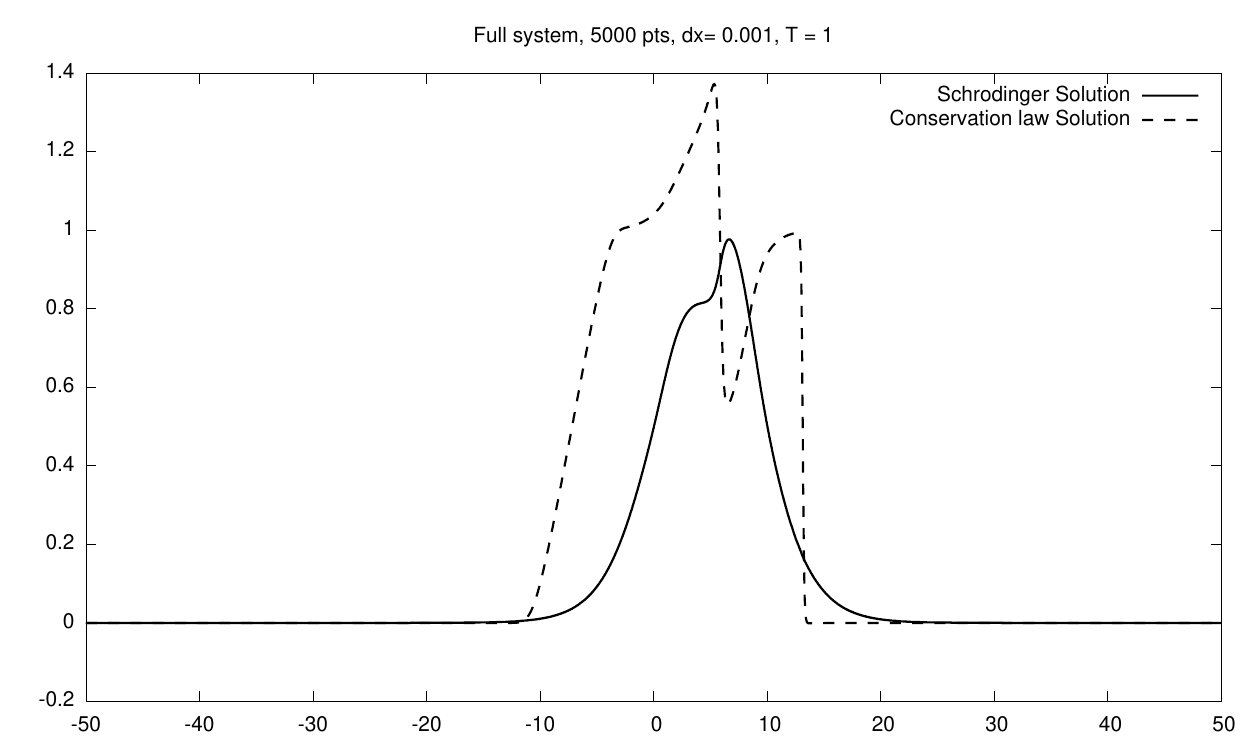}
\includegraphics[width=\linewidth,keepaspectratio=false]{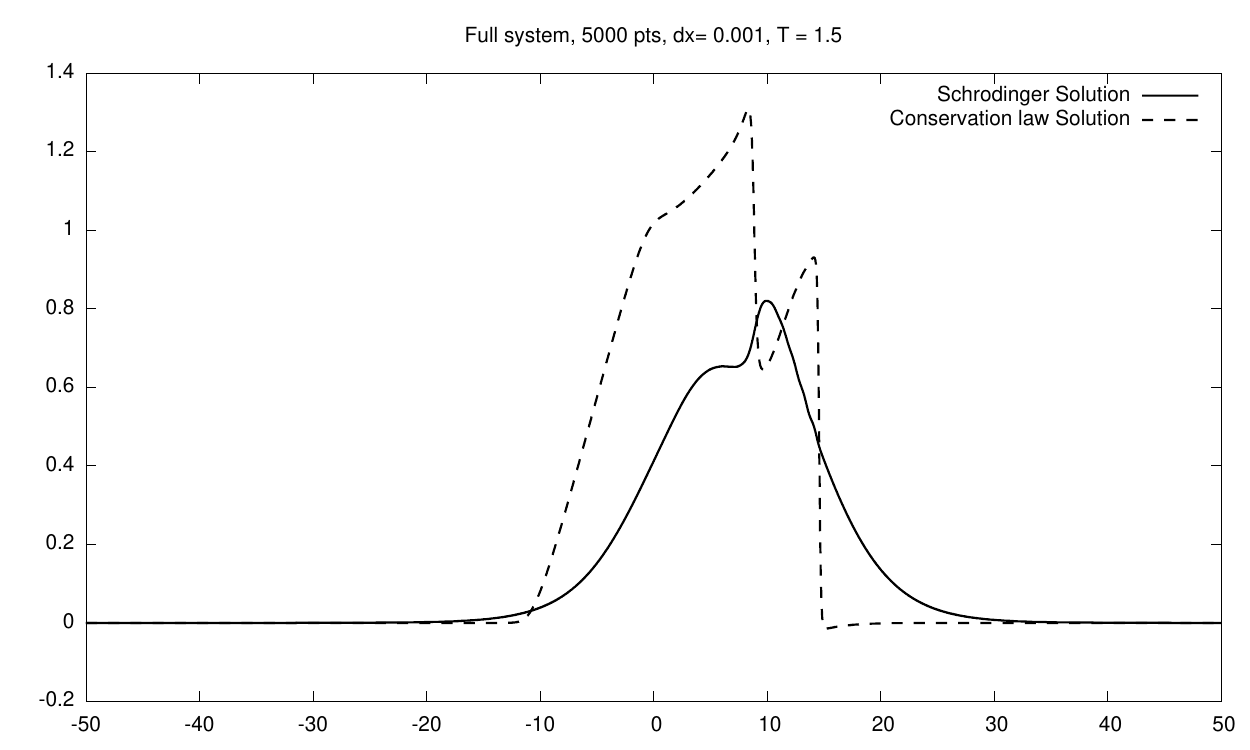}
\caption{Computed solutions of system \eqref{I.10},\eqref{I.20}, $T = 1,$ $T=1.5$.}
\label{FN-02}
\end{figure}

\begin{figure}
\includegraphics[width=\linewidth,keepaspectratio=false]{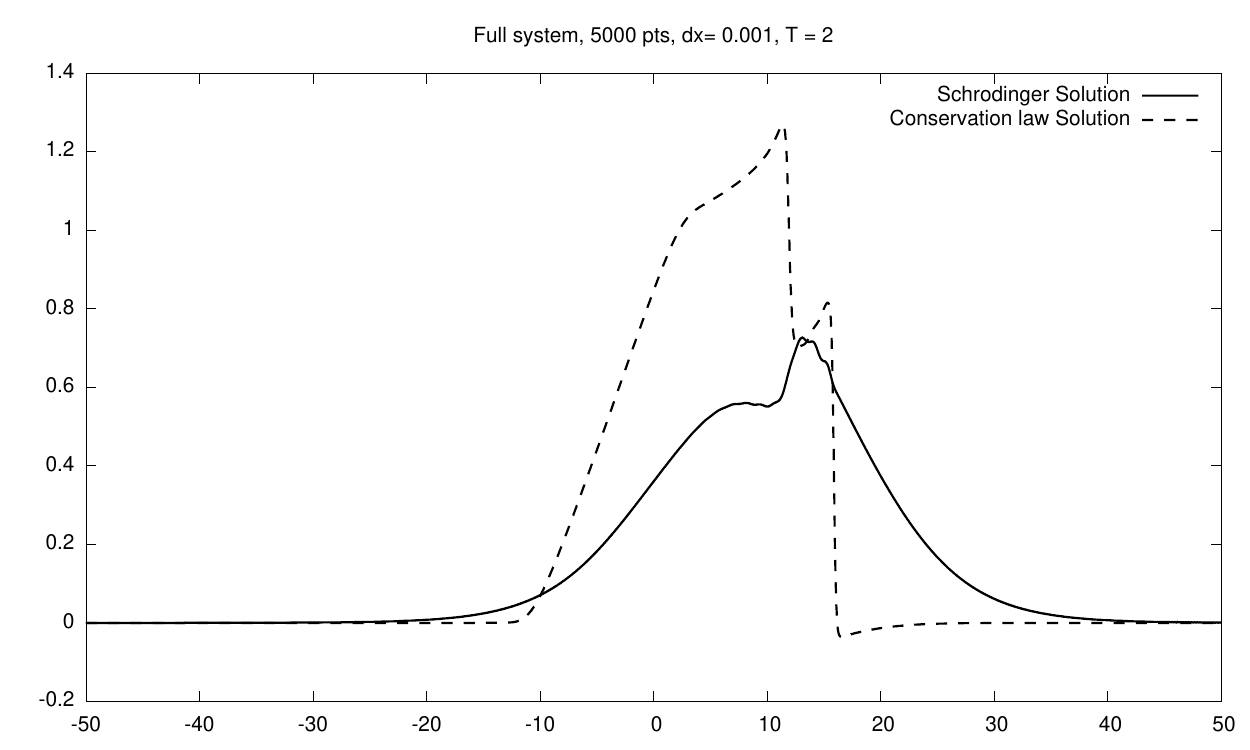}
\includegraphics[width=\linewidth,keepaspectratio=false]{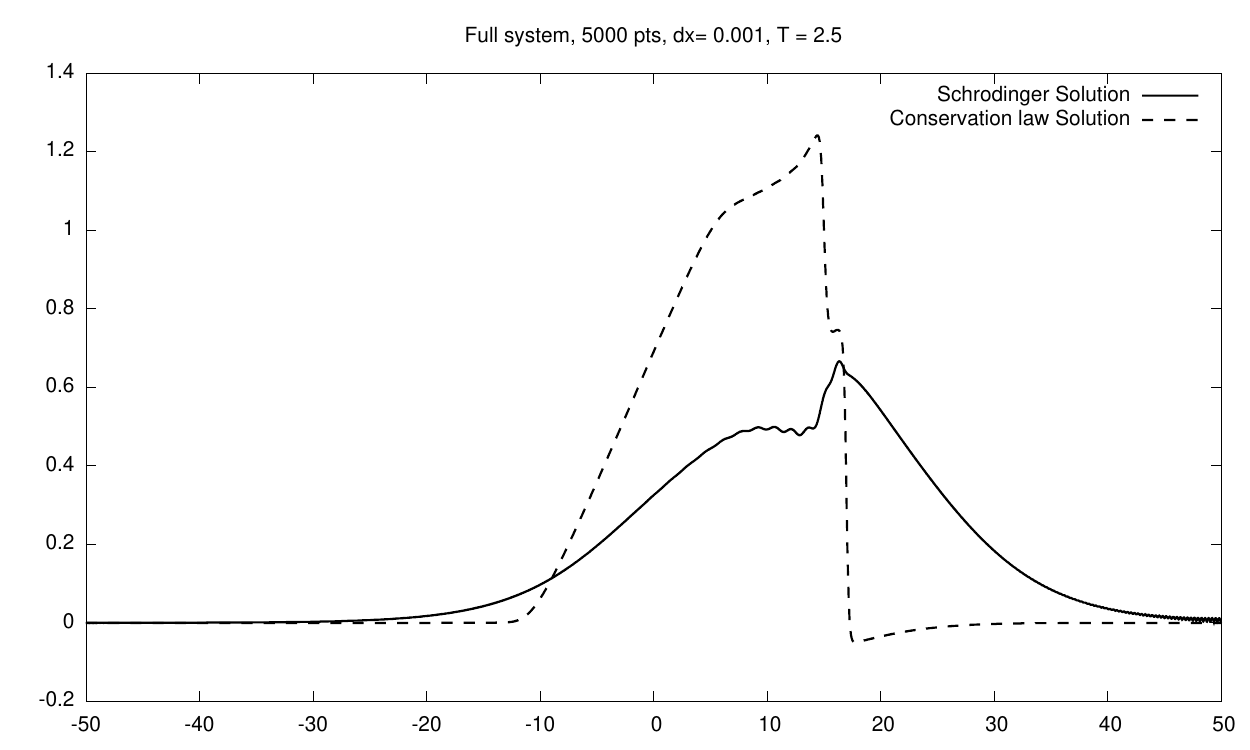}
\caption{Computed solutions of system \eqref{I.10},\eqref{I.20}, $T = 2,$ $T=2.5$.}
\label{FN-06}
\end{figure}

We take as initial data the functions
\[
\aligned
&u_0(x) = e^{5 ix/2}  \sqrt{6}\mathop{\mathrm{sech}}(\sqrt{3} x),
\\
&v_0(x) = \chi_{[-10,10]}
\endaligned
\]
on the spatial domain $[-50,50]$.
In Figures \ref{FN-01}--\ref{FN-06} we can observe
the usual behavior of the solutions to a nonlinear conservation law, such as the 
propagation of the initial shock discontinuity and the formation of rarefaction waves.
In addition, the interaction between the two equations induces the formation of new waves,
as can be seen clearly in Figure~\ref{FN-02}.

\begin{figure}
\includegraphics[width=\linewidth,keepaspectratio=false]{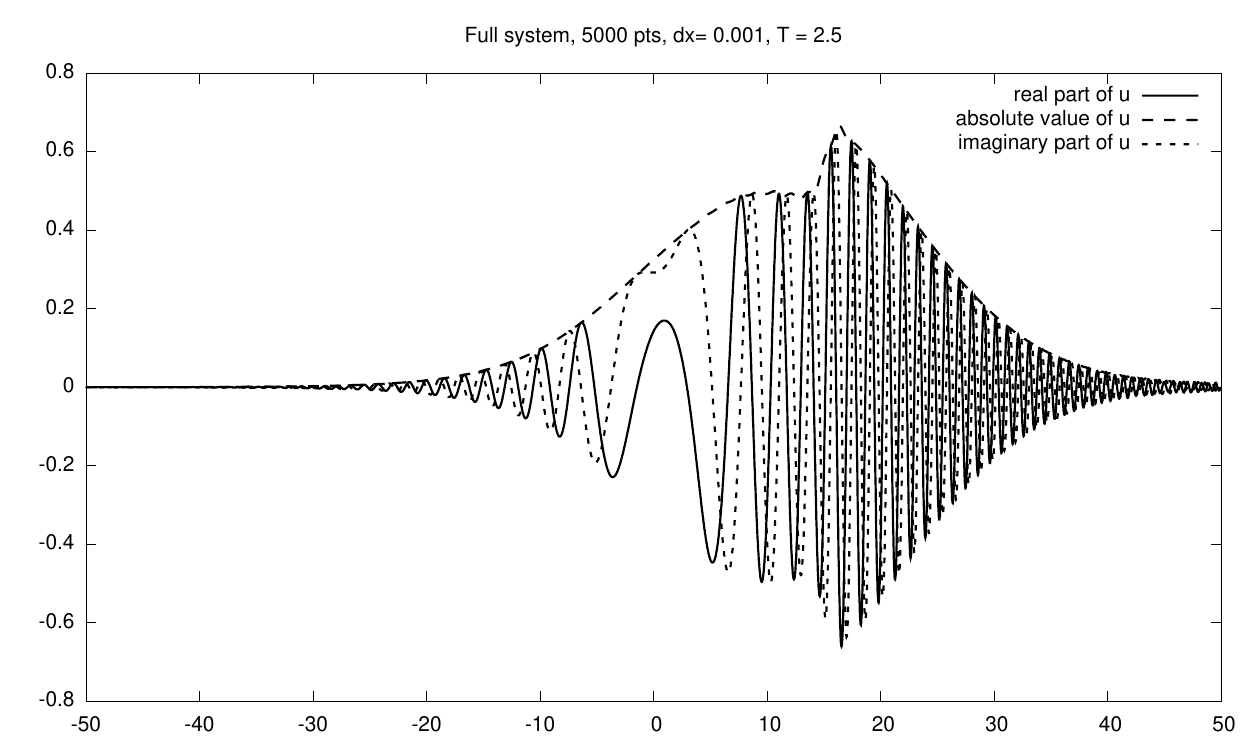}
\caption{Computed solution of system \eqref{I.10},\eqref{I.20}.}
\label{FN-07}
\end{figure}

For completeness, we show in Figure \ref{FN-07} the real and imaginary parts of the solution
$u$ of the Sch\-r\"o\-din\-ger equation, at time $T=2.5$.

Again, in these simulations, there was no need to consider the coupling function
$g(v)$, which is set to a cutoff function for some large cutoff parameter, since the value of
the solution $v$ never approaches this cutoff value.

Recall that one of the purposes of the coupling 
function $g$ was to ensure an $L^\infty$ bound on the solution $v$, which is in turn essential
to our convergence (and well-posedness) proof. Such an $L^\infty$ bound is an open question when
$g(v) = v$ \cite{ADFO,DF}. Now, in our experiments, we observed
no blow-up phenomena in $v$.
Since in our simulations the solution $u$ was always very regular, 
this suggests that, when the function $g$ is
the identity, any eventual blow-up of $v$ should be associated with initial data for $u$ (and thus
source terms for the conservation law \eqref{I.20}) having the minimum allowed regularity,
namely $u \in H^1(\RR)$ but $u_x \not\in L^\infty(\RR)$. A thorough study of this case,
however, is out of reach of our simple code,
and would need a more sophisticated numerical approach.

\section*{Acknowledgements}
PA would like to thank Nicolas Seguin and Nuno Lopes for useful discussions. The authors 
were partially supported by the FCT grant PTDC/MAT/110613/2009. PA was also supported by 
FCT through a {Ci\^encia 2008} fellowship. MF was also supported by FCT, 
Financiamento Base 2008-ISFL-1-297.


\end{document}